\documentclass[preprint,12pt]{elsarticle}

\usepackage{amssymb}
\usepackage{graphics}
\usepackage{graphicx}
\usepackage{subfigure}
\usepackage{epsfig}
\usepackage{multirow}
\usepackage[color,final]{showkeys}
\newtheorem{definition}{Definition}

\usepackage[aboveskip=1pt,labelfont=bf,labelsep=period,justification=raggedright,singlelinecheck=off]{caption}

\usepackage[normalem]{ulem}
\usepackage[color,final]{showkeys}
\usepackage[dvipsnames]{xcolor}
\usepackage{multirow}

\newcommand{\comm}[1]{{\color{red}(#1)}}

\journal{Information Sciences}
\begin{document}

\begin{frontmatter}



\title{Enhanced Joint Sparsity via Iterative Support Detection}


\author[label1]{Ya-Ru~Fan}
\author[label1,label2,label3]{Yilun~Wang \corref{cor1}}
\author[label1]{Ting-Zhu~Huang}
\address[label1]{ School of Mathematical Sciences, University of Electronic Science and Technology of
China, Chengdu, Sichuan, 611731 China}
\address[label2]{ PrinceTechs LLC., Shenzhen, Guangdong, 518101  P. R. China. }
\address[label3]{Center for Applied Mathematics, Cornell  University, Ithaca, NY, 14853, USA}
\cortext[cor1]{Corresponding author: yilun.wang@rice.edu (Yilun Wang) \\
E-mail addresses: yarufanfan@163.com (Ya-Ru Fan), yilun.wang@rice.edu (Yilun Wang), tingzhuhuang@126.com (Ting-Zhu Huang)}

\begin{abstract}
Joint sparsity has attracted considerable attention in recent years in many fields including sparse signal recovery in compressive sensing, statistics, and machine learning. Traditional convex models with joint sparsity suffer from the suboptimal performance though enjoying tractable computation.
In this paper, we propose a new non-convex joint sparsity model, and  develop a corresponding multi-stage adaptive convex relaxation algorithm.
This method  extends the idea of  iterative support detection (ISD) from the single vector estimation to the multi-vector estimation by considering the joint sparsity prior.
We provide some preliminary theoretical analysis including convergence analysis and a sufficient recovery condition.
Numerical experiments from both compressive sensing and multi-task feature learning show  the better performance of  the proposed  method in comparison with several state-of-the-art alternatives.
Moreover, we demonstrate that the extension of ISD from the single vector to multi-vector estimation is not trivial. While ISD doesn't well reconstruct the single channel sparse  Bernoulli signal, it does achieve significantly improved performance when recovering the multi-channel sparse  Bernoulli signal thanks to its ability of natural incorporation of the joint sparsity structure.

\end{abstract}


\end{frontmatter}


\section{Introduction and Contributions}
In the last decade, sparsity made it possible for us to reconstruct the high dimensional data with just few samples or measurements.  The key of the sparse estimation problem is to stress the identification of the support, which denotes the indices of the nonzeros. If the support is known, the estimation of the sparse vectors reduces to a standard overdetermined linear inverse problem \cite{Lee12Subspace}.

In order to enhance the estimation, many recent studies tend to consider the  structure information of the solutions. For example, 
group sparsity structure \cite{Group2011,GoncalvesDCSZB14}  widely appears in many applications \cite{singh2017block,zhu2016exploiting,7875091,sanders2017composite}, where the components of solutions are likely to be either all zero or all nonzero in a group.  By employing the grouping prior, ones  aim to  decrease the dispersion  to  facilitate  recovering a much better solution. Here, we focus on joint sparsity, which is a special case of the group sparsity. Joint sparsity means that  multiple unknown  sparse vectors share a common unknown nonzero support set. Unlike the many group sparsity situations where the grouping information is unknown, the joint sparsity provides us the group information.  In the following section, we will introduce the background of joint sparsity via two important applications, i.e. compressive sensing and multi-task feature learning.

In compressive sensing, joint sparsity aims to reconstruct unknown signals  from $m$ measurement vectors based on  a common measurement matrix. This is also called the multiple measurement vectors (MMV) problem \cite{Cotter05, Duarte05, Wimalajeewa14}. Given the observation vectors $b_{j}\in \mathbb{R}^{m}$ ($j=1,\cdots,\ell$) and a measurement matrix $A\in \mathbb{R}^{m\times n}$, we want to recovery the signal $x_{j}\in \mathbb{R}^{n}$ from the noisy underdetermined systems $b_{j}=Ax_{j}+e_{j}$, where $e_{j}\in \mathbb{R}^{m}$ is the noise. All the signal vectors $x_{1},\ldots,x_{\ell}$ share the sparsity pattern $M$, which implies the nonzero entries of $x_{1},\ldots,x_{\ell}$ almost appear on the same position. A common signal recovery model is
\begin{equation}\label{Model:L0}
\min_{x_{j}} ~~|M|   ~~~s.t. ~~ b_{j}=Ax_{j}+e_{j}, \quad j=1,\ldots,\ell,
\end{equation}
where $|M|$ is the cardinality of $M$ \cite{ReHB}.
In theory, we can recover the signals $X=[x_{1},\ldots,x_{\ell}]$ with rank $rank(X)=K$ if and only if
\begin{equation}
|M|<\frac{spark(A)-1+K}{2},
\end{equation}
where $spark(A)$ is the smallest set of linearly dependent columns of $A$ \cite{JATMJ06}. Since problem (\ref{Model:L0}) is NP-hard, it is usually relaxed with a convex alternative  which is computationally efficient at the cost of more required measurements. Like $\ell_1$-norm being the convex relaxation of $\ell_0$-norm \cite{Candes08IRL1}, the $\ell_{2,1}$-norm is widely used as the convex replacement of $|M|$ as below:
\begin{equation}\label{eq:l21_0}
\min_{X}~~~||X||_{2,1}:=\sum_{i=1}^n ||x^i||_{2}
\end{equation}
\begin{displaymath}
s.t.~~~b_j=Ax_j+e_j, j=1,\ldots, \ell,
\end{displaymath}
where $x^i \in \mathbb{R}^\ell$ and $x_{j} \in \mathbb{R}^n$ denote the $i$-th row and the $j$-th column of $X$, respectively.

Several fast algorithms have been proposed to solve problem (\ref{eq:l21_0}) \cite{7484756} such as greedy pursuit methods, iterative shrinkage algorithm \cite{Beck09} and alternating direction method (ADM) \cite{HLXL11}. The greedy pursuit methods such as matching pursuit and orthogonal matching pursuit (OMP) \cite{Wimalajeewa14} tend to require fewer computations but at the expense of slightly more measurements.

Multi-task learning has attracted much attention in machine learning \cite{fang2015visual,yang2017discriminative}. It aims to learn the shared information among related tasks in order for the improved performance than considering each learning task individually. Recently, multi-task feature learning based on the $\ell_{2,1}$-norm regularization has been studied. An underlying property of the $\ell_{2,1}$-norm regularization is that it urges multiple features from different tasks to share similar sparsity patterns \cite{Yang09Heterogenous}.
Given $\ell$ learning tasks associated with training data $\{(A^1,b_1)$, $\cdots$ , $(A^{\ell},b_{\ell})\}$, where
$A^{j} \in \mathbb{R}^{m_{j} \times n}$ is the data matrix of the $j$-th task with each row as a sample and each column as a feature; $b_j \in \mathbb{R}^{m_j}$ is the response  of the $j$-th task with biases $e_j$; $n$ is the number of features; and $m_j$ is the number of samples for the $j$-th task, we
would like to learn a weight matrix (sparsity pattern) $X = [x_1, \cdots,x_{\ell}] \in  \mathbb{R}^{n\times \ell}$ ($x_j \in \mathbb{R}^n $  consists of the
weight vectors for $\ell$ linear predictive models $b_j  = A^j x_j+e_j$ ) by solving the following optimization problem:
\begin{equation}\label{eq:l21_0-ml}
\min_{X}~~~||X||_{2,1}:=\sum_{i=1}^n ||x^i||_{2}
\end{equation}
\begin{displaymath}
s.t.~~~b_j =A^j x_j+e_j, j=1,\ldots, \ell,
\end{displaymath}
where $x^i \in \mathbb{R}^\ell$ and $x_{j} \in \mathbb{R}^n$ denote the $i$-th row and the $j$-th column of $X$, respectively. In this situation we assume that these different tasks share the same significant features, which leads to a joint sparsity problem.

The unconstrained formula corresponding to  problems (\ref{eq:l21_0}) and (\ref{eq:l21_0-ml}) can be written as the following unified form
\begin{equation}\label{eq:l21_unconstrained2}
\min_{X}~~~L(X) + \rho ||X||_{2,1},
\end{equation}
where $\rho>0$ is the regularization parameter, and $L(X)$ is a smooth convex loss function such as the least square loss function or the logistic loss function. For example, $L(X)= \sum^{\ell}_{j=1}\|A\textbf{x}_{j}-\textbf{b}_{j}\|_{2}^{2}$ for problem (\ref{eq:l21_0}), and $L(X)= \sum^{\ell}_{j=1}\frac{1}{\ell m_j}\|~A^j\textbf{x}_{j}-\textbf{b}_{j}\|_{2}^{2}$ for problem (\ref{eq:l21_0-ml}).

While the convexity of $\ell_{2,1}$-norm regularization provides computational efficiency, it also gives rise to  the inherited bias issue.  Similar with the $\ell_1$-norm regularized model 
which only achieves suboptimal recovery performance compared with the original cardinality regularized model from the theoretical viewpoint \cite{RCY08},   $\ell_{2,1}$-norm regularized model  also only achieves suboptimal performance compared with the cardinality based model  (\ref{Model:L0}).

Recently, several computational advances have been made  in the non-convex sparse regularization since its  performance is better than that of the convex sparse regularization \cite{fan2017cartoon,fan2016compressive}. For instance, for the single vector recovery, the non-convex $\ell_p$-norm ($0<p <1$) based sparsity regularization  usually obtains better performance than $l_{1}$-norm based sparsity regularization \cite{Charchtrand10lp,hu2016lp,jiang2016l_p}. For the joint sparsity, the  $\ell_{q,p}$-norm is applied in a similar way, where $0<p <1$ and $q \geq1$. The non-convex sparse regularization  needs less strict recovery requirements  and  usually achieves a  better  performance than the convex alternatives. While there have existed many algorithms for  solving  the non-convex sparse regularized models,  it is still a  challenging problem to  obtain  the  global  optimal  solution  efficiently.  The  behavior of a local solution is hard to analyze and more seriously  structural  information  of  the  solution is also hard to be incorporated into these algorithms.

$\mathbf{Contribution:}$ To achieve a better tradeoff between the recovery quality and the computational efficiency, we propose a non-convex joint sparsity regularized model and a multi-stage convex relaxation algorithm to solve the model. Motivated by the iterative support detection (ISD) \cite{Wang10ISD} for sparse signal reconstruction, we extend the idea of ISD in our method from common sparsity to joint sparsity, from compressive sensing to feature learning. We present some new insights about why ISD achieves better performance than its convex alternatives,  its key differences with other weighting based alternatives, and its flexibility in support detection implementation.  Moreover, we provide the preliminary theoretical results including the convergence analysis  and a sufficient recovery condition.

More importantly, we discover some advantages of ISD which are not observed in the single vector recovery. In particular, for the single channel sparse signal estimation, ISD depends on the assumption of the fast decaying property of the nonzero components of the underlying true sparse signal and does not work for non-decaying signals.  However, we empirically show that this assumption is no longer necessary for  multi-channel sparse signal recovery, because the joint sparsity structure is adopted in the specific implementation of support detection. This implies that ISD might be naturally fused with the general structural sparsity, which leads to the enhanced performance.

$\mathbf{Organization:}$ The remainder of this paper is organized as follows: in Section \ref{Sec:ours}, we propose a non-convex joint sparsity model and a corresponding algorithm based on ISD to solve the model. In Section \ref{sec:theory}, some preliminary theoretical results are presented. In Section \ref{Sec:NumExp}, we show numerical experiments on both compressive sensing and multi-task feature learning to demonstrate the effectiveness of the proposed method. Section \ref{sec:conclusion} is devoted to the conclusion and future works.

\section{The Proposed Model and Corresponding Algorithm} \label{Sec:ours}
We take the model (\ref{eq:l21_0}) of compressive sensing as the example to show how ISD is extended to the joint sparsity model. Similar idea can be expanded to the multi-task feature learning model (\ref{eq:l21_0-ml}), as well as the unconstrained version (\ref{eq:l21_unconstrained2}).

\subsection{Truncated Joint Sparsity Model}

The proposed  model based on the original joint sparsity model (\ref{eq:l21_0}) is given as follows:
\begin{equation} \label{truncated-Jointsparsity}
\min_{X(\omega)}~~~||X||_{w,2,1}:=\sum_{i=1}^n w_{i}||x^i||_{2}
\end{equation}
\begin{displaymath}
s.t.~~~B=AX+E,
\end{displaymath}
where $B=[b_{1},\ldots, b_{\ell}]$ is the observation matrix, $E$ is the noise matrix and $w = [w_{1}, w_{2},\cdots, w_{n}]$ is a weight parameter vector.  Compared with the model (\ref{eq:l21_0}), the main difference is the introduction of  the weight vector $w$.
Note that our model (\ref{truncated-Jointsparsity}) prefers a specific 0-1 weighting scheme, i.e.  $w_{i}$ is either $0$ or $1$, though most existing weighted models define weight as positive continuous real values.

Let $T$ be the set of the indices of the nonzero rows of $X$, and the model (\ref{truncated-Jointsparsity}) can be rewritten as
\begin{equation} \label{truncated-Jointsparsity_T}
\min_{X(T)}~~~||X||_{T,2,1}:=\sum_{i \in T} ||x^i||_{2}~~~~(TJS)
\end{equation}
\begin{displaymath}
s.t.~~~B=AX+E.
\end{displaymath}
We call it as truncated joint sparsity (TJS) model.

Intuitively, if we believe that $x^i$ is true nonzero,  it should be not forced to move closer to $0$ and therefore we  remove it from the regularization term, i.e. its corresponding $w_{i}$ is set as $0$.
While many existing works assume that  partial support information about underlying true sparse signal is already known \cite{Ince13,Saleh15}, the assumption may not hold in practice because $w_i$ is not given beforehand.  ISD,  as a self-learning scheme, aims to gradually detect partial support information. It is a multi-stage alternative optimization procedure, which repeatedly executes the following two steps when applied to model (\ref{truncated-Jointsparsity_T}):

 $\bullet$ Step 1: we optimize $x^{i}$ with $w$ (or $T$) fixed (initially $\vec{1}$): this is a convex problem in terms of $X$.

 $\bullet$ Step 2: we update $w$ using the current $X$ as reference via a support detection operation. The $w$ will be used in the Step $1$ of the next iteration. \\

Step 2 estimates the true nonzero rows from  the rough intermediate estimated results of Step 1, and therefore called ``support detection".
Our algorithm starts from initializing $w^{(0)}=\vec{1}$. In the first iteration, we obtain a solution $X^{(1)}$, which is the solution obtained by solving the plain $\ell_{2,1}$ model (\ref{eq:l21_0}). Then we achieve the weight $w^{(1)}$ using $X^{(1)}$ as the reference.
In the following iterations, we refine the intermediate solutions with the updated weights.
In fact ISD decouples the estimation of $w$ and $X$ by an alternative scheme. We denote this multi-stage convex relaxation procedure as iterative support detection based joint sparsity algorithm (ISDJS).

\subsection{Step 1: Solving Truncated Joint Sparsity  Model}

The $\ell_{2,1}$-norm based joint sparsity model (\ref{eq:l21_0}) leads to a convex optimization problem, and there are many  efficient first-order algorithms to solve it in different application fields \cite{Md09,KLY12,ZKDG13}, which  mostly  try to make use of the sparsity of the solutions in varied ways.
In compressive sensing, one of the most popular algorithms is the ADM method \cite{JYY09,Group2011}.
In \cite{JYY09}, Yong et al. used the ADM technique to solve the $\ell_{1}$-norm based optimization problem for compressed sensing and developed the corresponding Matlab package termed Your ALgorithms for $L_1$ (YALL1). Furthermore, Deng et al.  extended the YALL1 to the group version for solving the group sparse optimization with  $\ell_{2,1}$-norm regularization in \cite{Group2011}.
For feature learning, Liu et al. proposed an efficient algorithm based on the Nesterov¡¯s method and the Euclidean projection in \cite{JLiu}.

It is quite straightforward to extend these methods from  $\ell_{2,1}$-norm based models to truncated or weighted $\ell_{2,1}$-norm regularized models. We take the  YALL1 group algorithm for solving the plain $\ell_{2,1}$ regularized compressive sensing model as an example, and the resulted variant of the YALL1 group algorithm for the truncated joint sparsity  model (\ref{truncated-Jointsparsity}) is summarized in Algorithm 1,  where $\Lambda_{1}\in R^{n\times l}$, $\Lambda_{2}\in R^{m\times l}>0$ are multipliers in the ADM method, $\beta_{1}, \beta_{2}>0$ are penalty parameters, $\gamma_{1}$, $\gamma_{2}>0$ are step lengths, and $Z:=X$ is an auxiliary variable.

\begin{table}[!hb]
\begin{tabular}{l}
\hline
\textbf{Algorithm 1} Solving step 1 (inner loop)  \\
\hline
1.Initialize $X\in R^{n\times l}$, $\Lambda_{1}\in R^{n\times l}$, $\Lambda_{2}\in R^{m\times l}>0$,\\
  $\beta_{1}, \beta_{2}>0$ and $\gamma_{1}$, $\gamma_{2}>0$;\\
2.While stopping criterion is not met, do \\
(a)$X$ $\leftarrow$ $(\beta_{1}I + \beta_{2}A^TA)^{-1}(\beta_{1}Z -\Lambda_{1}+\beta_{2}A^{T}B + A^{T}\Lambda_{2})$,\\
(b)$Z\leftarrow$ Shrink $(X + \frac{1}{\beta_{1}}\Lambda_{1},\frac{1}{\beta_{1}}w)$ ,\\
(c)$\Lambda_{1}\leftarrow$ $\Lambda_{1} - \gamma_{1}\beta_{1}(Z - X)$ ,\\
(d)$\Lambda_{2}\leftarrow$ $\Lambda_{2} - \gamma_{2}\beta_{2}(AX - B)$,\\
where $\Lambda_{1}$, $\Lambda_{2}$ are multipliers, $\beta_{1}$, $\beta_{2}$ are penalty parameters, \\
      $\gamma_{1}$, $\gamma_{2}$ are step lengths.\\
\hline
\end{tabular}
\end{table}

The only  modification of the extension of the YALL1 group algorithm from the common joint sparsity to the truncated joint sparsity  is the step of updating $Z$, which is implemented by a shrinkage operator:
\begin{equation}
z^i=Shrink (r^i, \frac{1}{\beta_{1}}w)= max\{||r^i||_{2} - \frac{w_{i}}{\beta_{1}},0\}\frac{r^i}{||r^i||_{2}}, ~i=1,\cdots,n,
\end{equation}
where
\begin{equation}
r^i:=x^i + \frac{1}{\beta_{1}}\lambda_{1}^i.
\end{equation}

It is well known that $\ell_{2,1}$-norm based model, as its counterpart of $\ell_1$-norm based model, suffers from its uniform shrinkage on all its components, i.e., it shrinks the true nonzero components as well, and reduces the sharpness of the solution or introduces bias to the final solution.  In fact, the true  nonzero components should not be shrunk in order to avoid the possibly caused bias.   The truncated model can partially reduce this bias, because it corresponds to a selective shrinkage procedure where the weight value $w_i$ is either $1$ or $0$.  The true large nonzero components are expected to have the $0$ weights and thus will not be shrunk. Surely we need to have some knowledge about the  support information of the underlying true solution in order for the appropriate settings of weights. The support detection is implemented in Step 2, which will be introduced in the following  subsection. 

\subsection{Step 2: Weight Determination Based on the Iterative Support Detection} \label{sec:sd}
Step 2 is a vital part of the proposed algorithm. As mentioned above, our strategy obtains the partial support information by itself, rather than given beforehand.  Concretely, based on the recent intermediate result, we infer the indexes of nonzero rows of the underlying unknown true solution $\bar{X}$.
Once we believe that certain rows are nonzero in the true solution $\bar{X}$, we set the corresponding weights to be zeros, and the rest weights are all ones.

Since the intermediate results are not very accurate, a robust way to detect the correct information about the true nonzero rows is required. Some extra prior knowledge of the underlying $\bar{X}$ is  needed in order for reliable support detection. Recall that in the single channel sparse signal recovery case,  the nonzero components of the sparse or compressible signal  are assumed to have a fast decaying distribution for the effectiveness of the threshold based support detection scheme ($threshold$-ISD).  
As for the multi-vector estimation problem,  $threshold$-ISD, can also be applied in a similar way. At the $s$-th stage, we have an intermediate solution $X^{(s)}$. We aim to obtain some correct support information about the true $\bar{X}$ based on $X^{(s)}$, i.e. identify some truly nonzero rows.
The set of indices of detected nonzero rows  based on  $threshold$-ISD is similarly defined as follows:
\begin{equation}\label{eq:sd}
I^{(s+1)}:=\{i:|t_{i}^{(s)}|>\epsilon^{(s)}\}, s=0,1,2,\cdots,
\end{equation}
where $t_{i}=\|x^{i}\|_{2}$.  The support set $I^{(s)}$ is not necessarily increasing and nested, which means that all $s$ may be not in $I^{(s)}\subset I^{(s+1)}$. Because $I^{(s)}$ is not required to be monotonic, support detection can remove previous wrong detections, which makes $I^{(s)}$ less sensitive to $\epsilon^{(s)}$.

For the choice of $\epsilon^{(s)}$, the ``first significant jump'' heuristic was proposed in the original implementation of ISD \cite{Wang10ISD}.
Specifically, ones first sorts sequence $|t_{[i]}^{(s)}|$ in ascending order ($|t_{[i]}|$ denotes the $i$-th largest component of $t$ by magnitude). The ``first significant jump" scheme looks for the smallest $i$ such that
\begin{equation} \label{eq:support}
|t_{[i+1]}^{(s)}|-|t_{[i]}^{(s)}|>\tau^{(s)},
\end{equation}
where  $\tau^{(s)}$ is a data-dependent prescribed value to detect the big jump of this sequence. There are several simple and heuristic methods to define $\tau^{(s)}$. For example,
one can set
 \begin{equation} \label{eq:tau}
 \tau^{(s)} = m^{-1}||t^{(s)}||_{\infty},
 \end{equation}
where $m$ is the number of measurements.  Then we set $\epsilon^{(s)}=|t_{[i]}^{(s)}|$.
Intuitively, the ``first significant jump" scheme works well since the true nonzero entries of  ${t}^{(s)}$ are large in magnitude and small in number, while the false ones are large in number and small in magnitude. Therefore, the magnitudes of the true nonzero entries are spread out, while those of the false ones are clustered.

{ Recall that for sparse Bernoulli signal, where the nonzero components have exactly the same magnitude and do not have the fast decaying property, $threshold$-ISD fails to achieve a better performance than its convex alternative in the single channel recovery. Namely $threshold$-ISD works well only for the sparse signal whose nonzero components have the fast decaying magnitudes as presented in \cite{Wang10ISD}. However, for the joint sparsity situation, $threshold$-ISD naturally incorporates this extra joint sparsity structure in the implementation of support detection and succeeds achieving a better recovery quality, as experiments illustrated below in Section \ref{Sec:NumExp}. Here we give an intuitive explanation. 
Indeed, if we consider to each column  of $\bar{X}$ separately and perform threshold based support detection on each column individually, lack of fast decaying easily results in a significant number of wrong detections together with correct detections. This is the reason why for a single channel recovery, threshold based support detection will not help achieve a better recovery performance than the plain $\ell_1$ model when sparse Bernoulli signals are recovered.    However, for multiple sparse vectors which  share the same  sparsity structure,  the detected true non-zero positions of each individual vector  belong to the same subset (which contains all the true non-zero rows of the true solution $\bar{X}$ ) while the  detected false nonzero positions of each individual  vector might be quite different. Therefore, $|t_{i}^{(s)}|$ corresponding to the truly nonzero rows are much more likely to be significantly larger than those corresponding to the false nonzero rows.  Therefore, if we adopt the formula (\ref{eq:sd}) and (\ref{eq:support}) to perform support detection,  a relatively  accurate  support detection  can  be expected with an appropriate choice of threshold value. In other words,  we use the shared sparsity structure in the support detection procedure, i.e. the formula  (\ref{eq:sd}) and (\ref{eq:support}).
  }

We need to point out the difference of the 0-1 weighting scheme with a popular weighting method proposed in \cite{Candes08IRL1} where the weight is determined as follows:
$$
w_i^{(s)}=\frac{1}{|x^{(s)}_i|+\xi},
$$
where the choice of $\xi >0$ is a key for the performance of its corresponding algorithm.  If the $\xi$ is too small, then too much noisy information is taken into consideration. If the $\xi$ is too big, much of the information about the true nonzero elements is filtered out. An appropriate way might be to  gradually decrease $\xi$ from a large number to a small one as $s$ increases. However, the determination of $\xi$ is hard. While we know $\xi$ should be data-adaptive, a feasible practical scheme to determine $\xi$ is not easy to design. On the contrary, the Step 2 is much easier to obtain a data-adaptive scheme. The overall procedure of ISDJS is summarized in Algorithm 2.

In addition, an advantage of the proposed method is that the implementation of support detection is very flexible.
Besides the above heuristic  (\ref{eq:tau}), one can try other alternative ideas. For example, one dimensional edge detection methods such as those proposed in \cite{Zhang13slopedetection,Canny86} can also be adopted to detect the ``first significant jump", i.e. determine an appropriate $ \tau^{(\ell)}$ value for (\ref{eq:support}).  In the following numerical experiments, while the heuristic rule (\ref{eq:tau}) will be mostly adopted.   

\begin{table}[!ht]
\begin{tabular}{l}
\hline
\textbf{Algorithm 2} The ISDJS algorithm (outer loop)  \\
\hline
Input: measurement matrix $A$ and observation matrix $B$ ,\\
1.Set the iteration number $s\leftarrow 1$ and initialize $w^{(0)}=\overrightarrow{1}$ ;\\
2.While stopping condition is not met, do  \\
(a)$X^{(s)}$ $\leftarrow$ solve problem (\ref{truncated-Jointsparsity_T}) via Algorithm 1 for $w = w^{(s-1)}$;\\
(b)$w^{(s)}$ $\leftarrow$ $T^{(s)}$:= $(I^{(s)})^C =\{1,2\ldots,n\}$ $\backslash$ $I^{(s)}$;\\
 $\quad$ where $I^{(s)}$ $\leftarrow$ compute  approach (\ref{eq:sd}) for $X=X^{(s)}$;\\
(c)$s$$\leftarrow$ $s+1$.\\
\hline
\end{tabular}
\end{table}

The main computational cost of Algorithm 2 stems from computing the $X$ in Algorithm 1. We only need to compute the matrix inverse or do the matrix factorization once, whose complexity per iteration is $\mathcal{O}(mn)$ \cite{Group2011}.
Although ISDJS is a multi-stage procedure, 
the iteration number is small empirically, like around $4$ as the following numerical experiments presented in Section \ref{Sec:NumExp}. Thus, the total complexity of the proposed method is approximately $\mathcal{O}(mn)$.

\section{Preliminary Theoretical Analysis} \label{sec:theory}
Some preliminary theoretical analysis including convergence analysis and a sufficient recovery condition are presented, for the proposed TJS model and the ISDJS algorithm.

\subsection{Convergence Analysis}
We assume that $A \in \mathbb{R}^{m \times n} $ and $A^{j} \in \mathbb{R}^{m_{j} \times n}$ follow the continuous probability distribution. This convergence analysis applies to both the compressive sensing and the multi-task feature learning situations.  
%
%
In fact,  ISDJS only runs a few steps and $s$ is smaller than $5$ in general, and these steps can be considered to determine a proper threshold value $\epsilon$ for support detection (\ref{eq:sd}).  For simplicity of proof, when considering the convergence analysis where $s$ goes to infinity, we assume that  the threshold value $\epsilon^{(s)}$ used in support detection (\ref{eq:sd})  is fixed as $\bar{\epsilon}$ when $s>\bar{s}$ ($\bar{s}=5$, for example). This assumption is slightly biased from the truth, because the threshold value could keep changing as the iteration of ISDJS proceeds. However, it is also reasonable and acceptable to certain degree because in practice, ISDJS only runs a very limited number of steps, i.e. ISDJS will stop when $s$ is not very big. Moreover, even as $s$ goes to infinity, the threshold value will not change much empirically.

The main idea of the following proof refers to \cite{MSMTFL}. However,  unlike \cite{MSMTFL} where  a truncated $\ell_{1,1}$ model is considered, we consider a truncated $\ell_{2,1}$ model. A locally linear approximation is presented as a preparation for the following convergence proof.

First, for any $\epsilon>0$, we consider the following unstrained  weighted $\ell_{2,1}$ regularized model corresponding to model (\ref{eq:l21_unconstrained2}):
\begin{equation} \label{eq:unconstrainedTL21}
\min_{X} L(X)+\rho \sum_{i=1}^n w_{i}||\textbf{x}^i||_{2}
\end{equation}
where $L(X)$ is a quadratic cost function of $X$ and $\rho (>0)$ is a parameter. In terms of ISDJS algorithm we denote $w_{i}=T(\|\textbf{x}^{i}\|_{2}< \epsilon)$ and $T(\cdot)$ denotes the $\{0,1\}$-valued indicator function, which is consistent with $T=I^{C}$ in Algorithm 2.

The solution of problem (\ref{eq:unconstrainedTL21}) is equivalent to the solution of the following problem:
\begin{equation} \label{eq:cappedL21}
\min_{X} L(X)+ \rho \sum_{i=1}^{n} \min(\|\textbf{x}^{i}\|_{2}, \epsilon)
\end{equation}

We assume that $\epsilon^{(s)}$  is kept unchanged after several stages of ISDJS. That is to say when $s$ is big enough,  $\epsilon^{(s)} \doteq \bar{\epsilon}$ in Algorithm 2.

Second, we define two auxiliary functions:
\begin{equation}
\textbf{h}: R^{n\times l}\longmapsto R^{n}_{+}, \textbf{h}(X)=[\|\textbf{x}^{1}\|_{2}, \|\textbf{x}^{2}\|_{2}, \cdots, \|\textbf{x}^{n}\|_{2}]^{T},
\end{equation}

\begin{equation}
g_{\epsilon}: R^{n}_{+}\longmapsto R_{+}, g_{\epsilon}(\textbf{u})= \sum_{i=1}^{n} \min(u_{i}, \epsilon).\qquad\qquad\qquad
\end{equation}
Note that $g_{\epsilon}(\cdot)$ is a concave function \cite{MSMTFL}. (When  $s$ is large enough, $ g_{\epsilon^{(s)}}=g_{\bar{\epsilon}}$.) We know that a vector $\textbf{z}\in R^{n}$ is a sub-gradient of $g$ at $\textbf{v}\in R_{+}^{n}$, if for all vector $\textbf{u}\in R_{+}^{n}$, the following inequality holds:
\begin{equation}\label{neq:1}
g_{\epsilon}(\textbf{u})\leq g_{\epsilon}(\textbf{v})+\langle \textbf{z}, \textbf{u}-\textbf{v}\rangle,
\end{equation}
where $\langle \cdot, \cdot \rangle$ denotes the inner product. Based on the functions defined above, problem (\ref{eq:cappedL21}) is equivalent to the following problem
\begin{equation}
\min_{X} L(X)+\rho g_{\epsilon}(\textbf{h}(X)).
\end{equation}
Then we can obtain an upper bound of $g_{\epsilon}(\textbf{h}(X))$ using a locally linear approximation at $\textbf{h}(X^{(s)})$ based on the inequality (\ref{neq:1}):
\begin{equation}
g_{\epsilon}(\textbf{h}(X))\leq g_{\epsilon^{(s)}}(\textbf{h}(X^{(s)}))+ \langle \textbf{z}^{(s)}, \textbf{h}(X)-\textbf{h}(X^{(s)})\rangle ,
\end{equation}
where $\textbf{z}^{(s)}= [T(\|(\textbf{x}^{(s)})^{1}\|_{2}< \epsilon), \cdots, T(\|(\textbf{x}^{(s)})^{n}\|_{2}< \epsilon)]^{T}$ is a sub-gradient of $g_{\epsilon}(\textbf{u})$ at $\textbf{u}= \textbf{h}(X^{(s)})$. 
Furthermore, for $\forall~\textbf{h}(X)$ we obtain an upper bound of the optimization problem (\ref{eq:cappedL21}):
\begin{equation} \label{neq:2}
L(X)+\rho g_{\epsilon}(\textbf{h}(X))\leq L(X)+\rho g_{\epsilon^{(s)}}(\textbf{h}(X^{(s)}))+\rho \langle \textbf{z}^{(s)}, \textbf{h}(X)-\textbf{h}(X^{(s)})\rangle.
\end{equation}
Since $\rho$ and $h(X^{(s)})$ are constant with respect to $X$, we have
$$X^{(s+1)}= \arg\min_{X}  L(X)+\rho g_{\epsilon^{(s)}}(\textbf{h}(X^{(s)}))+\rho \langle \textbf{z}^{(s)}, \textbf{h}(X)-\textbf{h}(X^{(s)})\rangle$$
\begin{equation} \label{neq:3}
 = \arg\min_{X}  L(X)+\rho (\textbf{z}^{(s)})^{T}\textbf{h}(X),\qquad\qquad\quad\qquad\quad
\end{equation}
which corresponds to the step 2 (a) of the Algorithm 2. Therefore, we intuitively consider that the Algorithm 2 minimizes an upper bound in each step.

\textbf{Theorem 1.}\emph{ Let  $\bar{f}(X)=  L(X)+\rho \sum_{i=1}^n \bar{w}_{i}||\textbf{x}^i||_{2}$, where $\bar{w}_{i}=T(\|\textbf{x}^{i}\|_{2}< \bar{\epsilon})$.  The sequence $\{\bar{f}(X^{(s)})\}$ is decreasing and convergent.}
\vspace{2mm}

\textbf{Proof.} When $s$ is big enough, $\epsilon^{(s)}=\bar{\epsilon}$, and $g_{\epsilon^{(s)}}=g_{\bar{\epsilon}}$. 
We firstly verify the objective function value in (\ref{eq:cappedL21}) decreases monotonically based on locally linear approximation, when $s$ is big enough, as follows:
$$L(X^{(s+1)})+\rho g_{\bar{\epsilon}} (\textbf{h}(X^{(s+1)}))\leq L(X^{(s+1)})+\rho g_{\bar{\epsilon}}(\textbf{h}(X^{(s)}))+\rho \langle \textbf{z}^{(s)}, \textbf{h}(X^{(s+1)})-\textbf{h}(X^{(s)})\rangle$$
$$\leq L(X^{(s)})+\rho g_{\bar{\epsilon}}(\textbf{h}(X^{(s)}))+\rho \langle \textbf{z}^{(s)}, \textbf{h}(X^{(s)})-\textbf{h}(X^{(s)})\rangle$$
$$= L(X^{(s)})+\rho g_{\bar{\epsilon}}(\textbf{h}(X^{(s)})),\qquad\qquad\qquad\qquad\qquad\qquad$$
where the first inequality is based on equation (\ref{neq:2}) and the second inequality follows (\ref{neq:3}), i.e.,  $X^{(s+1)}$ is a minimizer of the right hand side of equation (\ref{neq:2}).
Then we have
$$\bar{f}(X^{(s+1)})\leq \bar{f}(X^{(s)}).$$

Obviously, we observe that
$$\bar{f}(X)= L(X)+\rho \sum_{i=1}^n \bar{w}_{i}||\textbf{x}^i||_{2} \geq 0.$$
Thus $\{\bar{f}(X^{(s)})\}$ is bounded below.  Therefore $\{\bar{f}(X^{(s)})\}$ is convergent. $\blacksquare$ 

\subsection{A Sufficient Recovery Condition of TJS Model}
Here we discuss the noiseless  compressive sensing model (\ref{eq:l21_0}).
We first review the truncated null space property (t-NSP) proposed in \cite{Wang10ISD}, which is a generalization of the null space property (NSP).
\begin{definition} Matrix  $A$ satisfies the $t$-NSP of order $L$ for $\gamma>0$ and
$0<t \le n$ if
\begin{equation} \label{eq:TNSP}
\|\eta_{S}\|_{1} \le \gamma \|\eta_{(T \cap S^C)}\|_{1}
\end{equation}
holds for all sets $T\subset\{1,\ldots,n\}$ with $|T|=t$, all subsets $S
\subset T$ with $|S|\le L$, and all $\eta \in \mathcal{N}(A)$ --- the null space of $A$.

For simplicity, we use $t$-NSP$(t, L, {\gamma})$ to denote the $t$-NSP of  order $L$ for ${\gamma}$ and $t$, and use $\bar{\gamma}$ to replace $\gamma$ and write $t$-NSP$(t, L, \bar{\gamma})$ if $\bar \gamma$ is the infimum of
all the feasible $\gamma$ satisfying (\ref{eq:TNSP}).
\end{definition}

For the single channel sparse signal recovery problem \begin{equation}\label{eq:l21_0-single}
\min_{x}~~~||x_T||_{1}
\end{equation}
\begin{displaymath}
s.t.~~~b=Ax
\end{displaymath}
Theorem 3.2 in \cite{Wang10ISD} has shown that if $A\in \mathbb{R}^{m\times n} (m<n)$ satisfies the t-NSP, then the true signal $\bar{x}$ is the unique solution of model (\ref{eq:l21_0-single}) if
\begin{equation} \label{eq:kd}
||\bar{x}_{T}||_{0}<k(d),
\end{equation}
where $k(d):= c\frac{m-d}{1+log(\frac{n-d}{m-d})}$, $d=n-t=n-|T|$, and $c>0$ is absolute constant independent of the dimensions $m$, $n$ and $d$.
Let $d_{c}= |I\cap supp(\bar{t})|$ stand for the number of correct detections, and  inequality (\ref{eq:kd}) is equivalent to
\begin{equation} \label{eq:kd2}
||\bar{x}||_{0}< k(d)+d_{c},
\end{equation}
due to $||\bar{x}||_{0}= ||\bar{x}_{T}||_{0}+d_{c}$.
In light of (\ref{eq:kd2}), to compare the common $\ell_1$ model with truncated $\ell_1$ model (\ref{eq:l21_0-single}), we shall compare $k(0)$ with $k(d) + d_c$.
In \cite{Wang10ISD}, the authors have shown that if there are enough correct detections (i.e., $d_c/d$ is sufficiently large), then we get $k(0) < k(d)+d_c$.  That is to say, the truncated $\ell_1$ model might be able to recovery more nonzeros than the common $\ell_1$ model. 
We extend the above conclusion about the advantage of the truncated $\ell_1$ model over the common $\ell_1$ model from the single vector case to the multi-vector recovery case, i.e. the joint sparsity case. We will see that this kind of extension is feasible thanks to the theorem proved in \cite{NSP}, which is revisited as below.

Theorem 1.3 in \cite{NSP} has proved the following two statements are equivalent.

(a) for all vectors $(\textbf{u}^{1},\cdots, \textbf{u}^{r})\in ( \mathcal{N}(A))^{r}\setminus \{(0, 0, \cdots, 0)\}$
\begin{equation}
\sum_{j\in S}(\sqrt{u_{1,j}^{2}+\cdots+ u_{r,j}^{2}})< \sum_{j\in S^{c}}(\sqrt{u_{1,j}^{2}+\cdots+ u_{r,j}^{2}}),
\end{equation}

(b) for all vector $\textbf{v}\in  \mathcal{N}(A)$ with $v\neq 0$
\begin{equation}
\sum_{j\in S} |v_{j}|< \sum_{j\in S^{c}} |v_{j}|, \textbf{v}=(v_{1},\ldots, v_{n})^{T}\in R^{n}.
\end{equation}
They hold for all index sets $S\subset \{1,2,\ldots, n\}$ with $|S|\leq L$, where $ \mathcal{N}(A)$ stands for the null space of $A$ and $S^{c}$ is  the complement set of $S$. Namely, the null space property of multiple systems of linear equations is equivalent to the null space property for the comm $\ell_1$ minimization subject to a single linear system.

During their proof of this equivalence, they only make use of the fact $S \cap S^{c} =\emptyset$. So we naturally have the following equivalence \cite{NSP}:

(c) for all vectors $(\textbf{u}^{1},\cdots, \textbf{u}^{r})\in ( \mathcal{N}(A))^{r}\setminus \{(0, 0, \cdots, 0)\}$
\begin{equation}
\sum_{j\in S}(\sqrt{u_{1,j}^{2}+\cdots+ u_{r,j}^{2}})< \sum_{j\in T\cap S^{c}}(\sqrt{u_{1,j}^{2}+\cdots+ u_{r,j}^{2}}),
\end{equation}

(d) for all vector $\textbf{v}\in  \mathcal{N}(A)$ with $v\neq 0$
\begin{equation}
\sum_{j\in S} |v_{j}|< \sum_{j\in T\cap S^{c}} |v_{j}|, \textbf{v}=(v_{1},\ldots, v_{n})^{T}\in R^{n}.
\end{equation}
They hold for all index sets $S\subset \{1,2,\ldots, n\}$ with $|S|\leq L$. Thus, we similarly have the equivalence of the t-NSP of multiple systems of linear equations with the t-NSP  for the common $\ell_1$ minimization subject to a single linear system .

Therefore, the better recovery performance of the truncated $\ell_1$ model over the common $\ell_1$ model can be extended to  our specific joint sparsity case, i.e. the multiple-vector compressive sensing problem. In other words, if there are enough correct detections, the truncated joint sparsity model (\ref{truncated-Jointsparsity}) can recover more nonzero rows than the common $\ell_{2,1}$ model, based on the above simple and intuitive analysis. However, for multi-task learning problem, where different $A_j$ exists, the theoretical judgement of the truncated $\ell_1$ model over the common $\ell_1$ model needs further investigation.



\section{Numerical Experiments} \label{Sec:NumExp}
We show some numerical experiments to demonstrate the better performance of the proposed ISDJS in comparison with several state-of-the-art algorithms. For compressive sensing, the YALL1 group \cite{Group2011}, SOMP \cite{JATMJ06} and p-threshold \cite{Gribonval07Thresholding} algorithms are compared. For multi-task feature learning, ISDJS is compared with the baseline algorithm for the common $\ell_{2,1}$ regularized model, whose  baseline algorithm  proposed in \cite{JLiu} is implemented in the software ``MALSAR" \cite{Zhou12MALSAR}.  We mainly focus on the  recovery rate and accuracy. Due to the channel number has a great influence on recovery rate of joint sparsity, we provide several different channel number settings. In addition, we also test the robustness of the competing approaches in different noise levels.  The synthetic experiments and two realistic experiments on collaborative spectrum sensing \cite{J2011} and multi-task feature learning, verify the effectiveness of ISDJS.

\subsection{Parameter Settings of  Synthetic Examples}
Two synthetic examples are presented for compressive sensing. The true jointly sparse solution $\bar{X}\in R^{n\times \ell}$ is generated by randomly selecting $k$ rows as nonzero rows whose entries follow the i.i.d  Gaussian and Bernoulli distribution in test 1 and test 2, respectively. The rest rows of $\bar{X}$ are set as zeros. Randomized partial Walsh-Hadamard transform matrix is utilized as the measurement matrix $A\in R^{m\times n}$ in compressive sensing, because it is suitable for large-scale computation and has the property $AA^{T}=I$. For SOMP and p-thresholding algorithms, we use their default parameter settings in  \cite{JATMJ06,Gribonval07Thresholding}. We set the parameters of the YALL1 group algorithm and ISDJS referring to \cite{Group2011} as follows: $\beta_{1} = 0.3/\frac{1}{m\ell}\sum_{i=1}^{m}\sum_{j=1}^{\ell}|b_{ij}|$, $\beta_{2} = 3/\frac{1}{m\ell}\sum_{i=1}^{m}\sum_{j=1}^{\ell}|b_{ij}|$ ($b_{ij}$ is the entries of matrix $B\in R^{m\times \ell}$ ) and $\gamma_{1} = \gamma_{2} = 1.618$.
All involved algorithms are  terminated when
\begin{equation}
\frac{||t^{(k+1)}-t^{(k)}||_{2}}{||t^{(k+1)}||_{2}}\leq \varepsilon.
\end{equation}
For SOMP, p-threshold and YALL1 group algorithms, $\varepsilon$ is set as $10^{-6}$. As for ISDJS, in the first few outer loops, we only want to get an rough estimate of the support information of $X$, thus we just set a relatively loose tolerance such as $\varepsilon=10^{-2}$. But in the last iterations,   $\varepsilon$ is also set as $10^{-6}$ for fair comparison. In all experiments, ISDJS runs no more than $5$ outer loops.

The empirical recovery performance of all test algorithms  in general becomes  better as the number of channels gradually increases, though to varying degrees \cite{Eldar}. Therefore, different channel number settings are tried in our experiments, for example,  $L=1,2,4,8,16$, respectively. We also try different sparsity levels varying from $k=80$ to $160$, while fixing  $n=1024$, $m=256$ in all tests excluding Figs \ref{fig:1}, \ref{fig:2}, \ref{fig:6}, \ref{fig:7}. The experimental results corresponding to compressive sensing are usually an average of $100$ runs due to the involved randomness in the generation of $A$ and $\bar{X}$.

\subsubsection{Test 1: Compressive recovery of joint sparse Gaussian signals}

\begin{figure}[!htpb]
  \centering
   \subfigure[]{
   \includegraphics[width=6.6cm,height=3.6cm]{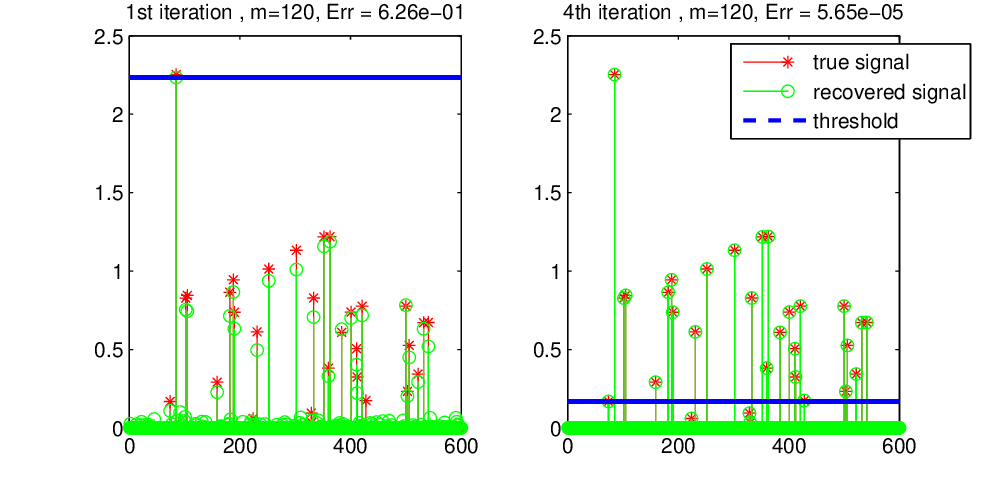}}
     \subfigure[]{
   \includegraphics[width=6.6cm,height=3.6cm]{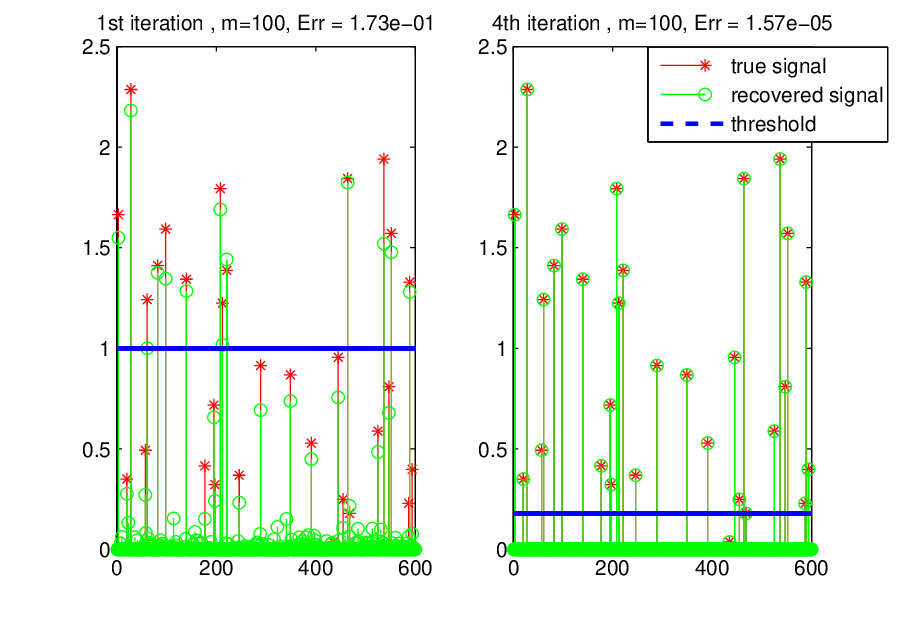}}
     \subfigure[]{
   \includegraphics[width=6.6cm,height=3.6cm]{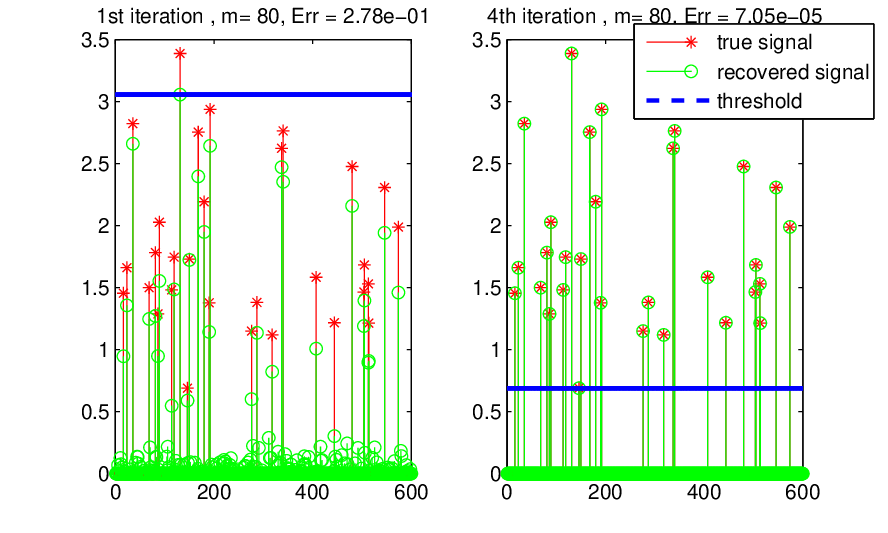}}
     \subfigure[]{
   \includegraphics[width=6.6cm,height=3.6cm]{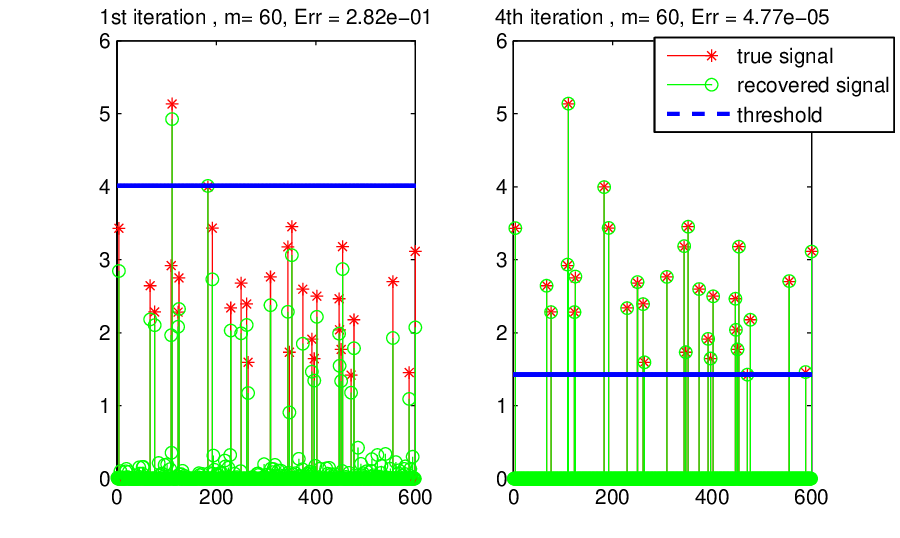}}
     \subfigure[]{
   \includegraphics[width=6.6cm,height=3.6cm]{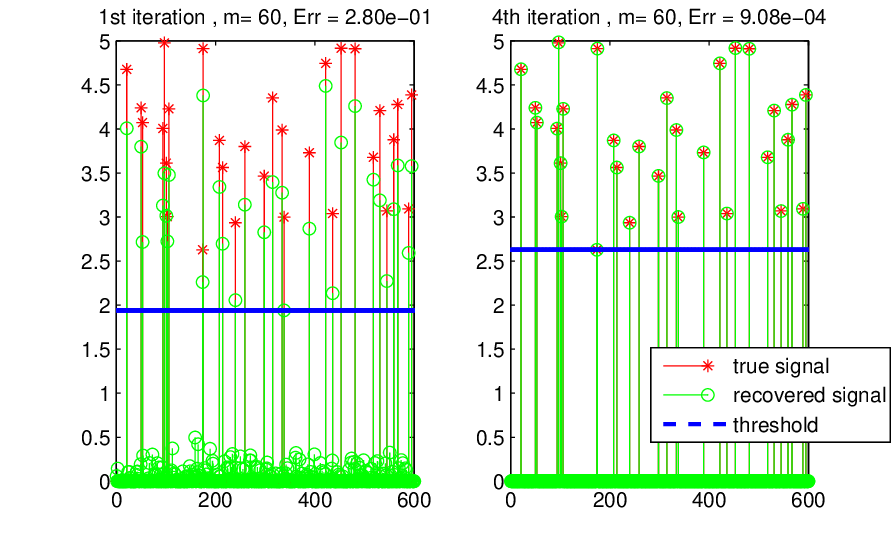}}
   \centering
\begin{tabular}{|c|c|c|c|c|c|}
\hline 
 \multirow{2}{*}{L--Iteration} & \multicolumn{4}{|c|}{Nonzeros} &  \multirow{2}{*}{Relative error}     \\
 \cline{2-5}
     & Total true & Detected & Correct & False &  \\
    \cline{1-6}
     1--1 &  30 &37 	& 29   &  8  & 6.26e-01 \\

     1--4 &  30 & 30	& 30   &  0  & 5.65e-05 \\
     \cline{1-6}

     2--1 &  30 & 34	& 28  &  6  & 3.73e-01 \\

     2--4 &  30 & 30 & 30  &  0  & 7.57e-05 \\
    \cline{1-6}
     4--1 &  30 & 38	& 30  &  8  & 2.78e-01 \\

     4--4 &  30 & 30 & 30  &  0  & 7.05e-05 \\
    \cline{1-6}
     8--1 &  30 & 34	& 30  &  4  & 2.82e-01 \\

     8--4 &  30 & 30 & 30  &  0  & 4.77e-05 \\
    \cline{1-6}
     16--1 &  30 & 33	& 30  &  3  & 2.80e-01 \\

     16--4 &  30 & 30 & 30  &  0  & 9.08e-05 \\
    \cline{1-6}
\hline
\end{tabular}
   \caption{\small Compare the true Gaussian signals and recovered signals from ISDJS in different channels, where the two parts in each subplot are the results in the first iteration and the fourth iteration respectively. (a)L=1, m=120, (b)L=2, m=100, (c)L=4, m=80, (d)L=8, m=60, (e)L=16, m=60.}
   \label{fig:1}
\end{figure}

We perform a synthetic compressive sensing example to demonstrate some key aspects of ISDJS including the effectiveness of threshold based support detection, the effect of the channel number. We also illustrate how ISDJS produces gradually improved intermediate solutions starting from the low quality initial point which is the solution of  the convex alternative. We generate a sparse signal $\bar{X} \in R^{600\times L}$ with $k=30$ nonzero rows. The results of ISDJS in the first iteration and the fourth iteration for different channel numbers settings are depicted in the Fig \ref{fig:1}, where we set $\bar{t}$ (a vector of 2-norm of each row of $\bar{X}$) on behalf of the true signal $\bar{X}$ and $t$ (a vector of 2-norm of each row of $X$ from ISDJS) on behalf of the recovered signal.  We use the quadruplet``(Total, Detected, Correct, False)" and ``Err"  defined below to  measure the accuracy of support detection.

$ \bullet$ (Total, Detected, Correct, False): \\
$ \quad -$ Total: the number of total nonzero rows of the true signal $\bar{X}$; \\
$ \quad -$ Detected: the number of detected nonzero rows, equal to $|I| = (Correct) + (False)$; \\
$ \quad -$ Correct: the number of correctly detected nonzero rows, i.e., $|I \cap \{i: \bar{t}^{i} \neq 0\}|$; \\
$ \quad -$ False: the number of falsely detected nonzero rows , i.e., $|I \cap \{i: \bar{t}^{i} = 0\}|$.

$ \bullet$ Err: the relative error $\|X-\bar{X}\|_{2}/\|\bar{X}\|_{2}$. \\

From Fig \ref{fig:1} (a), we can see that the output of the first iteration of ISDJS, which is the solution of the common convex $\ell_{2,1}$-norm regularized model (\ref{eq:l21_0}), is not good. Nevertheless the output of the fourth iteration of ISDJS could well match the true signal with a very small relative error.
We also see that ISDJS is insensitive to a small number of false detections and has an attractive self-correction capacity. In particular, while it is difficult for the common $\ell_{2,1}$ model (\ref{eq:l21_0}) to recover a signal with $k=30$ nonzero entries from $m=60$ measurements, ISDJS can finally return a satisfying result, as presented in Fig \ref{fig:1} (e). Note that when the measurements $m$ decrease, ISDJS still returns a better result with channel numbers $L$  increasing.

\begin{figure}[!htbp]
  \centering
    \subfigure[]{
   \includegraphics[scale=0.45]{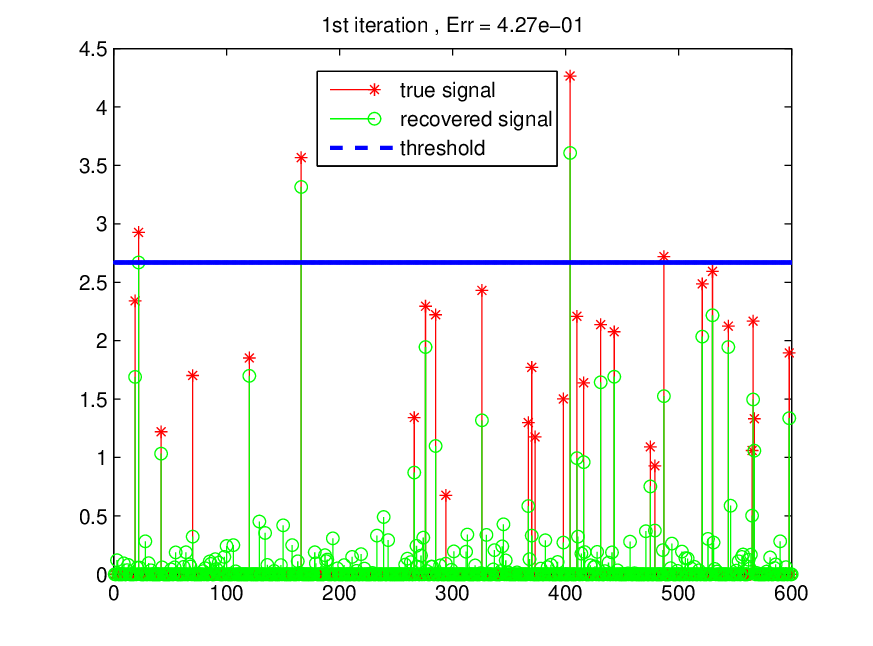}}
     \subfigure[]{
   \includegraphics[scale=0.45]{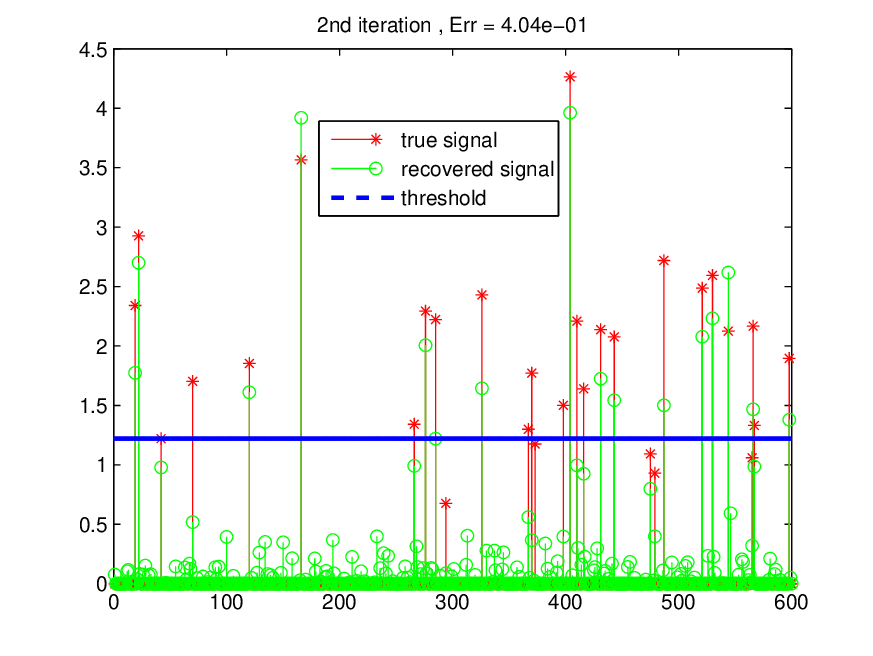}}
     \subfigure[]{
   \includegraphics[scale=0.45]{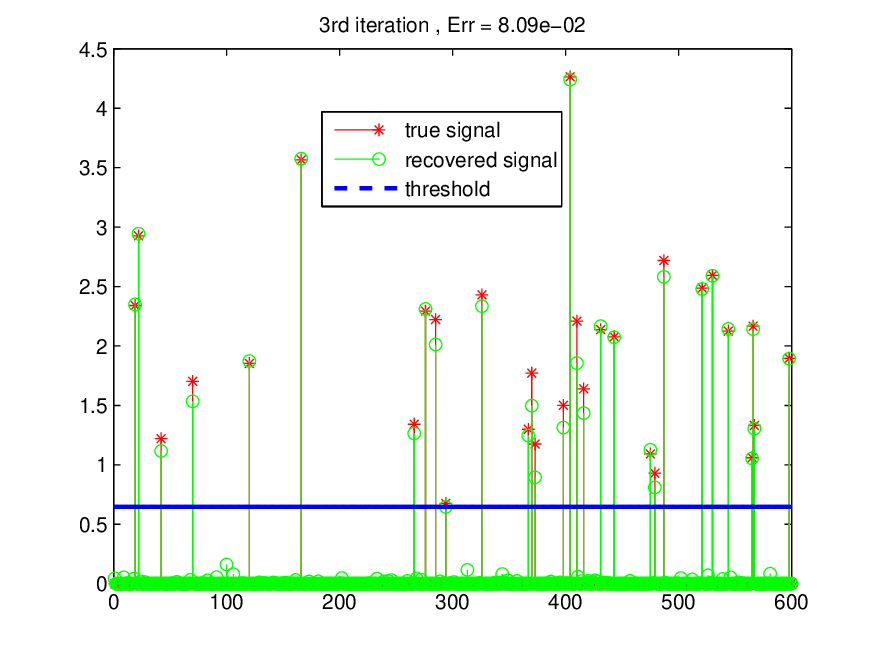}}
     \subfigure[]{
   \includegraphics[scale=0.45]{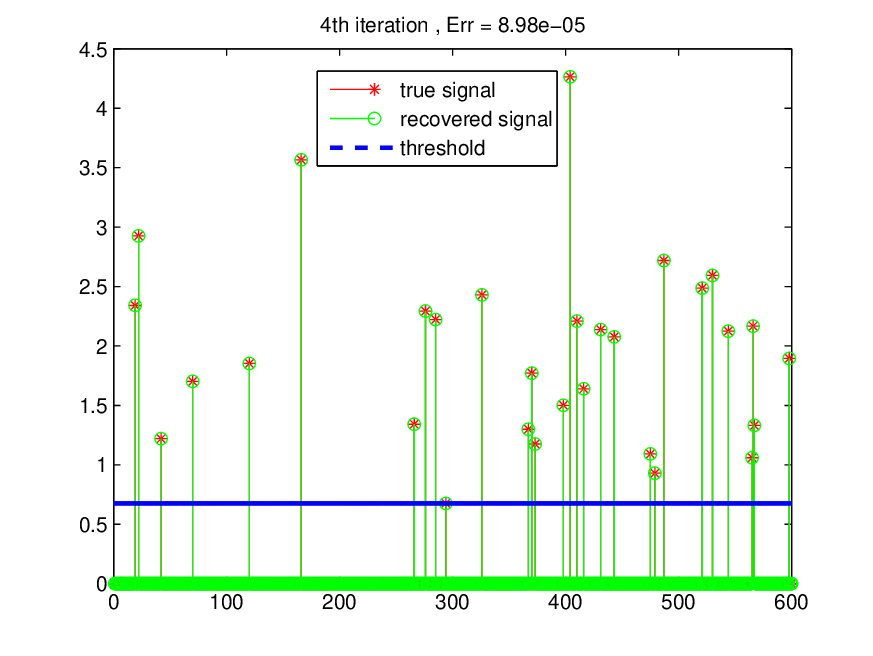}}
\centering
\begin{tabular}{|c|c|c|c|c|c|}
\hline 
 \multirow{2}{*}{Iteration} & \multicolumn{4}{|c|}{Nonzeros} &  \multirow{2}{*}{Relative error}     \\
 \cline{2-5}
     & Total true & Detected & Correct & False &  \\
    \cline{1-6}
     1 &  30 & 41	& 27   &  14  & 4.27e-01 \\

     2 &  30 & 40	& 28  &  12  & 4.04e-01 \\

     3 &  30 & 30	& 30  &  0  & 8.09e-02 \\

     4 &  30 & 30   & 30  &  0  & 8.98e-05 \\
    \cline{1-6}
\hline
\end{tabular}
   \caption{\small Compare the true Gaussian signals and recovered signals obtained by ISDJS in each iteration with L=4 channels.}
   \label{fig:2}
\end{figure}

In order to better understand ISDJS, we show each outer iteration of it in Fig \ref{fig:2},  by taking an example of $L=4$ and $m=80$. From the Fig \ref{fig:2} (a), ISDJS in the first iteration (i.e. YALL1 group algorithm), finds very few positions of correct nonzero rows and has a large relative error. However, a half  positions of correct nonzero rows could be detected in the next iteration as exhibited in Fig \ref{fig:2} (b), and  a significantly improved recovery is obtained, shown in Fig \ref{fig:2} (c).  In the third iteration, our algorithm  has already correctly detected the most nonzero positions, and therefore a good enough solution is obtained as illustrated in Fig \ref{fig:2} (d). 

\begin{figure}[!htbp]
  \centering
   \subfigure[]{
   \includegraphics[scale=0.44]{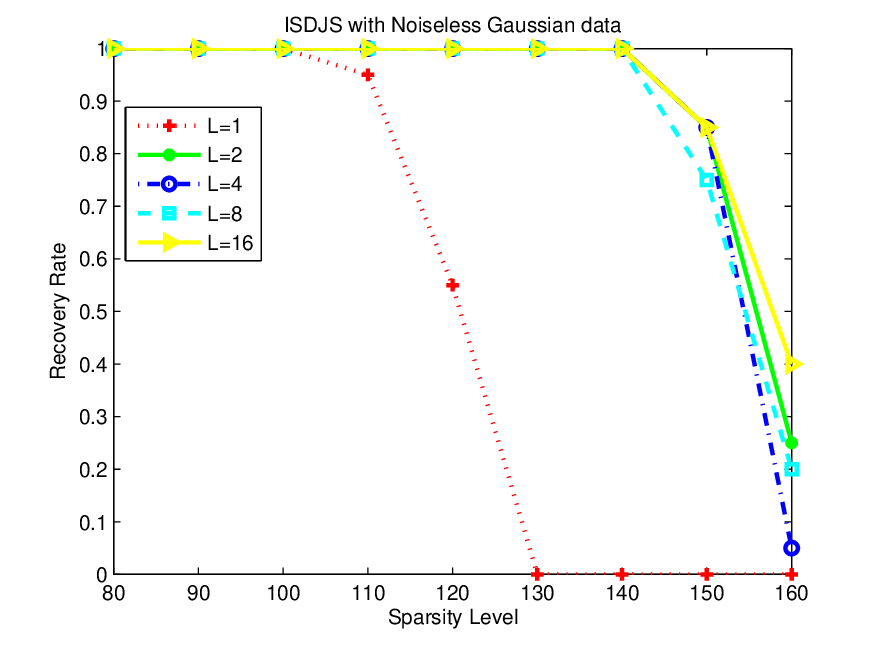}}
   \subfigure[]{
   \includegraphics[scale=0.44]{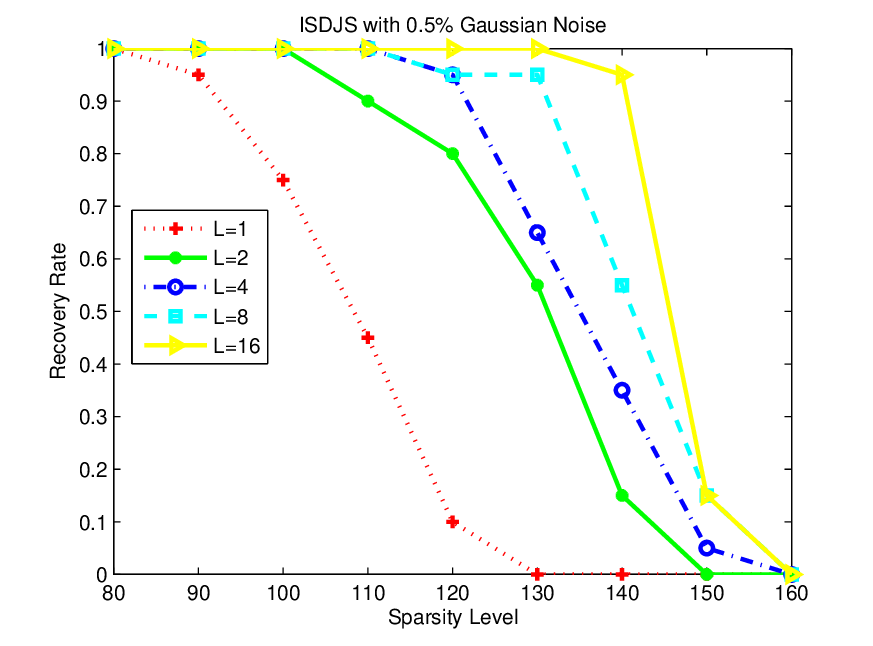}}
   \subfigure[]{
   \includegraphics[scale=0.44]{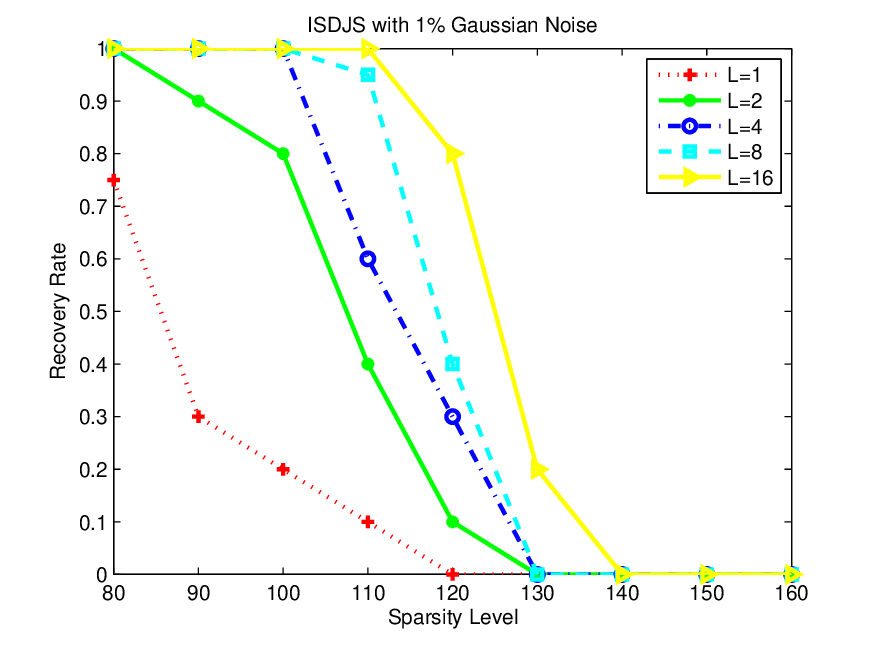}}
  \subfigure[]{
   \includegraphics[scale=0.44]{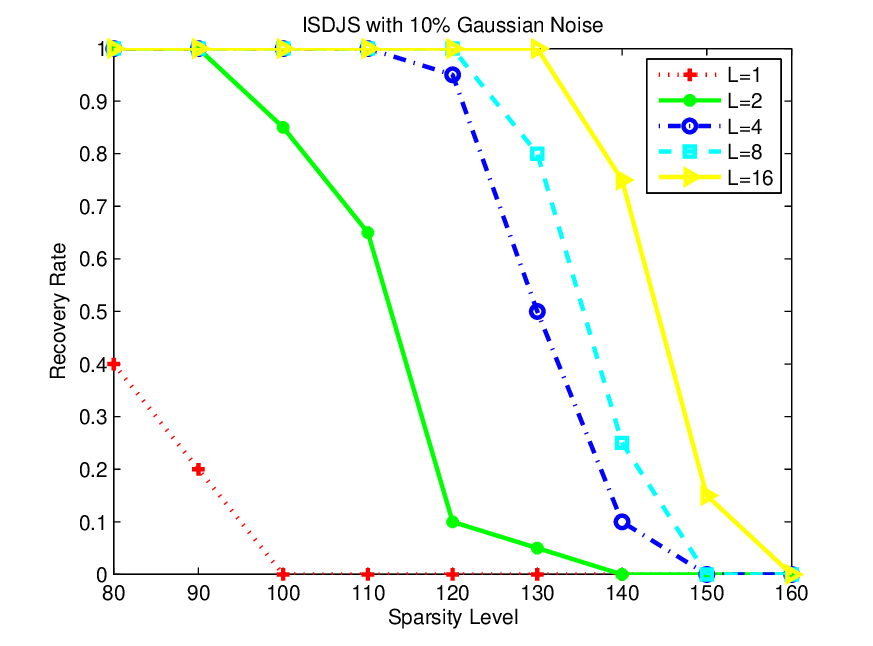}}
   \caption{\small Compare the recovery rate of ISDJS with L=1, 2, 4, 8, 16 in different noise levels for Gaussian signals, (a)noiseless, (b)0.5\% noise, (c)1\% noise, (d)10\% noise.}
   \label{fig:3}
\end{figure}

Fig \ref{fig:3} shows the performance of ISDJS with different channels in four different noise levels to verify its robustness. The proposed algorithm performs better for the multi-channel sparse signal recovery than the single channel sparse signal recovery even in the high level noise.

\begin{figure}[!htbp]
  \centering
   \subfigure[]{
   \includegraphics[scale=0.43]{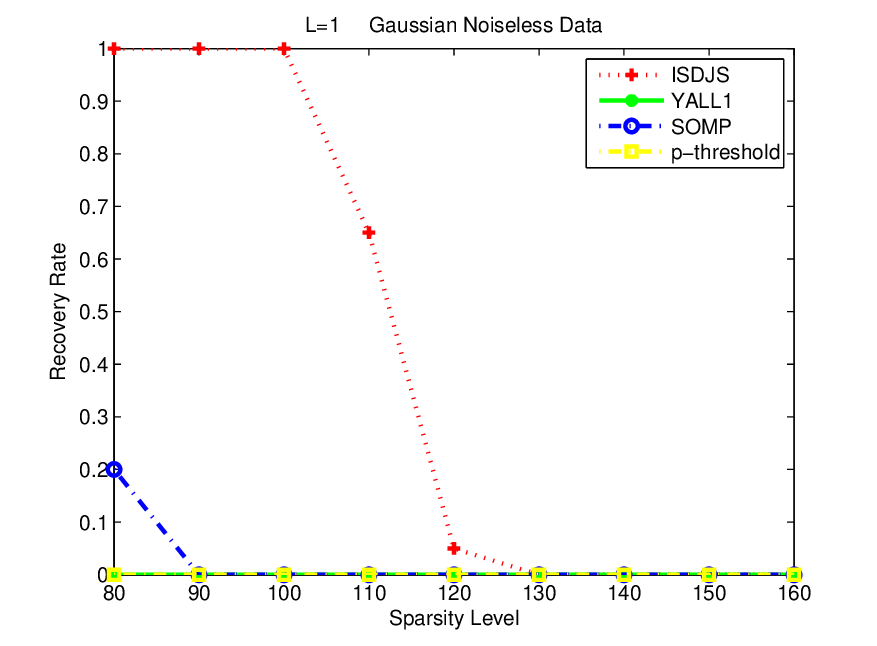}}
   \subfigure[]{
   \includegraphics[scale=0.43]{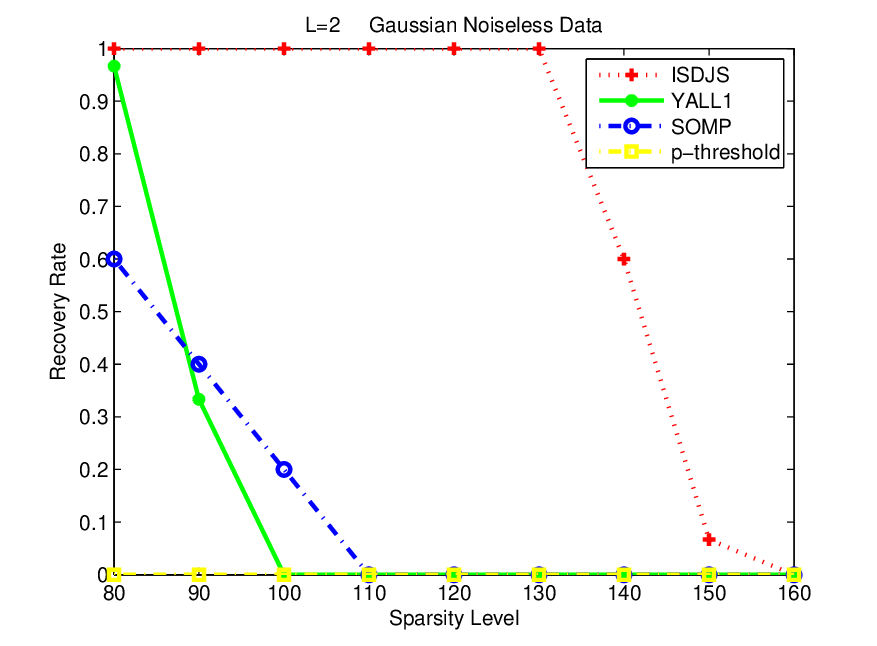}}
  \subfigure[]{
   \includegraphics[scale=0.43]{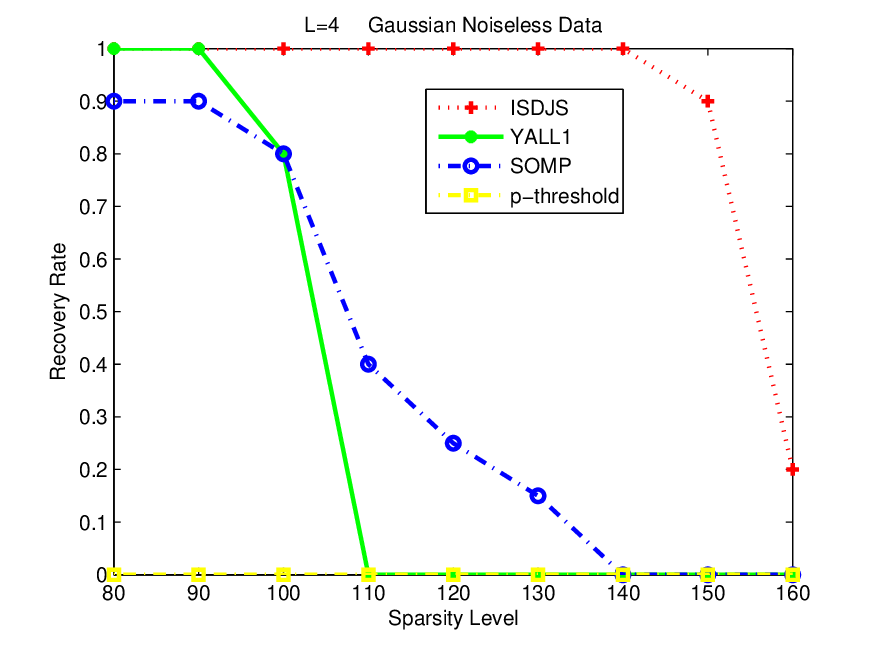}}
   \subfigure[]{
   \includegraphics[scale=0.43]{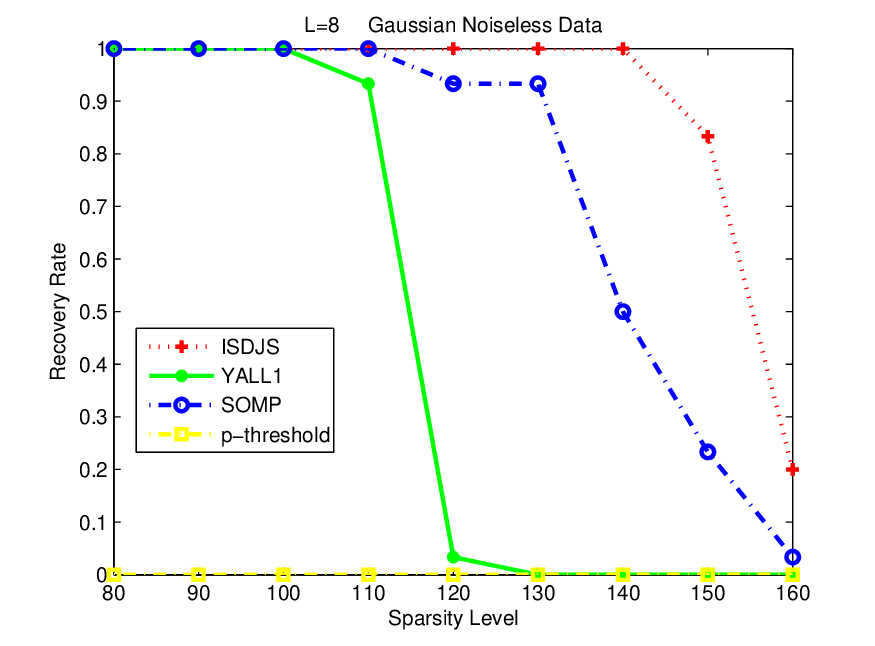}}
   \caption{\small Compare the recovery rate of four algorithms in different channels for noiseless Gaussian signals, (a)L=1, (b)L=2, (c)L=4, (d)L=8.}
   \label{fig:4}
\end{figure}

\begin{figure}[!htbp]
  \centering
   \subfigure[]{
   \includegraphics[scale=0.44]{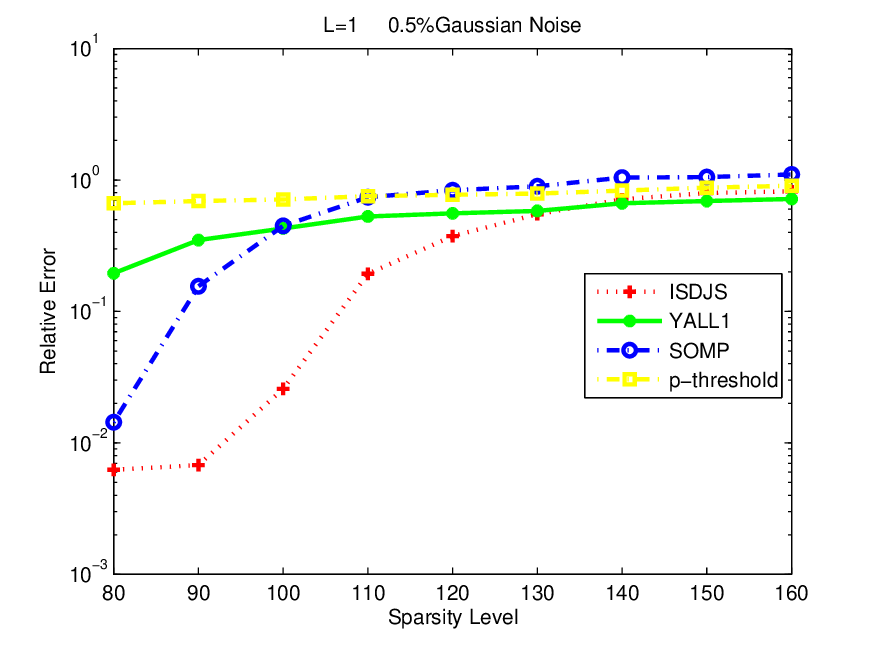}}
   \subfigure[]{
   \includegraphics[scale=0.44]{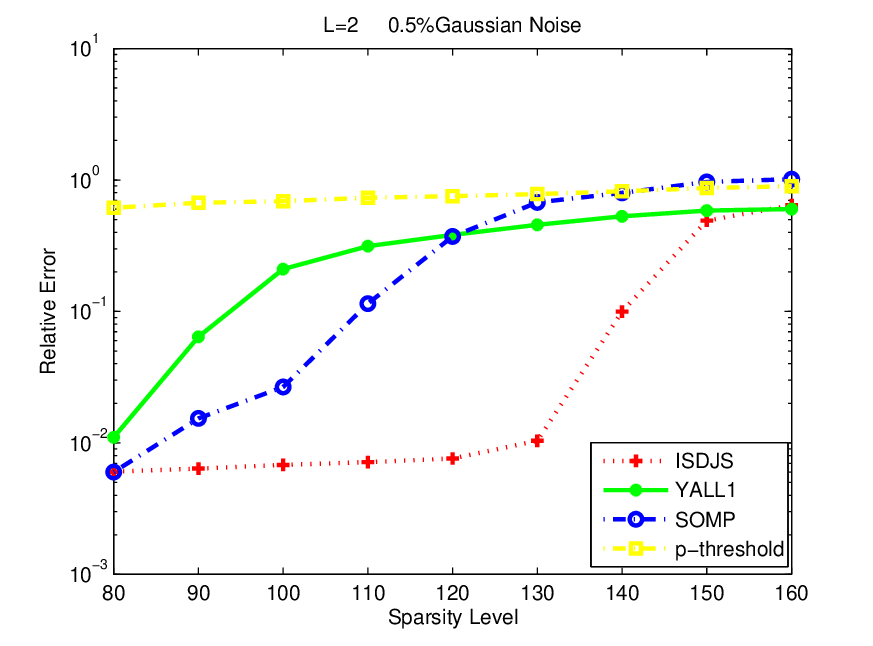}}
   \subfigure[]{
   \includegraphics[scale=0.44]{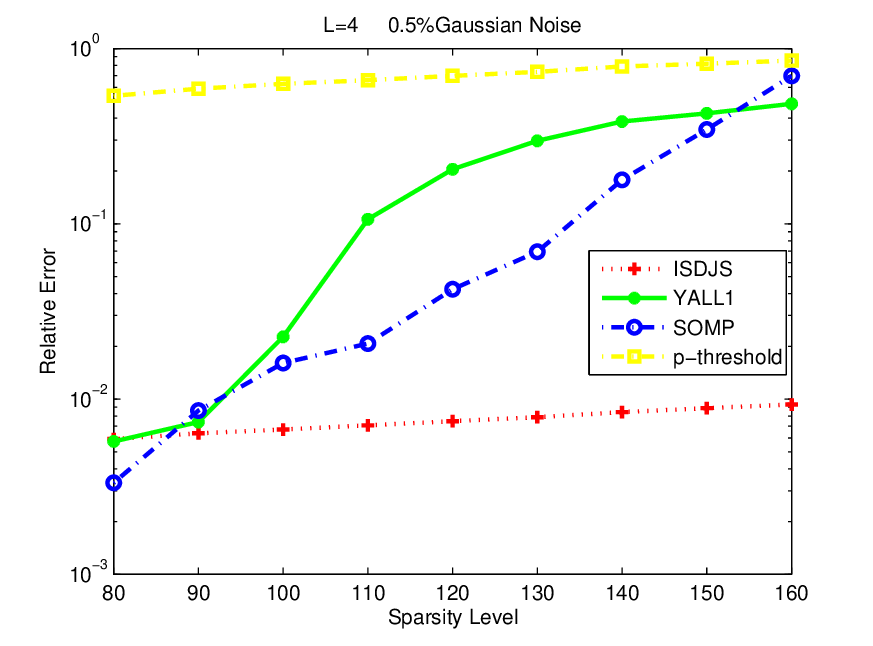}}
   \subfigure[]{
   \includegraphics[scale=0.44]{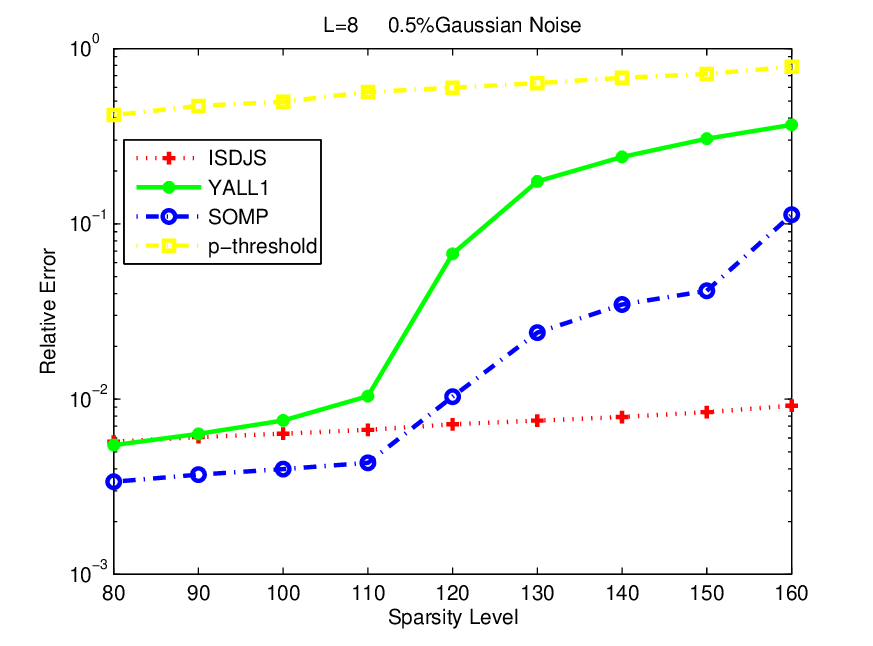}}
   \caption{\small Compare relative error of four algorithms in different channels with 0.5\% noise for Gaussian signals, (a)L=1, (b)L=2, (c)L=4, (d)L=8.}
   \label{fig:5}
\end{figure}

In Fig \ref{fig:4} and Fig \ref{fig:5}, we compare the recovery rates and relative errors of all test algorithms, for noiseless case and noisy case (added Gaussian noise with standard variance $0.5\%$), respectively, when the channel number varies.  We can see that ISDJS outperforms other algorithms in all involved different channels. While the common $\ell_{2,1}$ model behaviors worse than the SOMP in the cases of $L=4$ and $L=8$,  ISDJS which applies ISD to the common $\ell_{2,1}$ model, works better than SOMP.

\subsubsection{Test 2: Compressive recovery of joint sparse Bernoulli signals}
\begin{figure}[!htbp]
  \centering
   \subfigure[]{
   \includegraphics[scale=0.43]{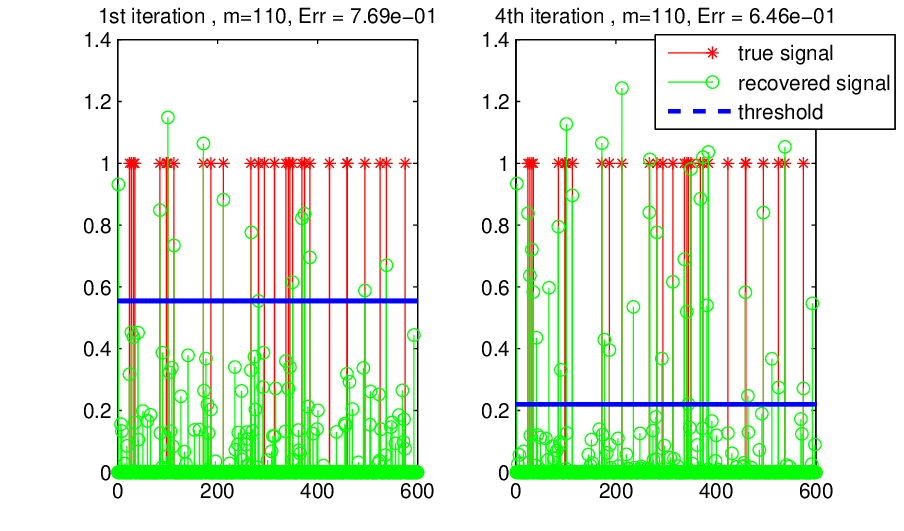}}
   \subfigure[]{
   \includegraphics[scale=0.43]{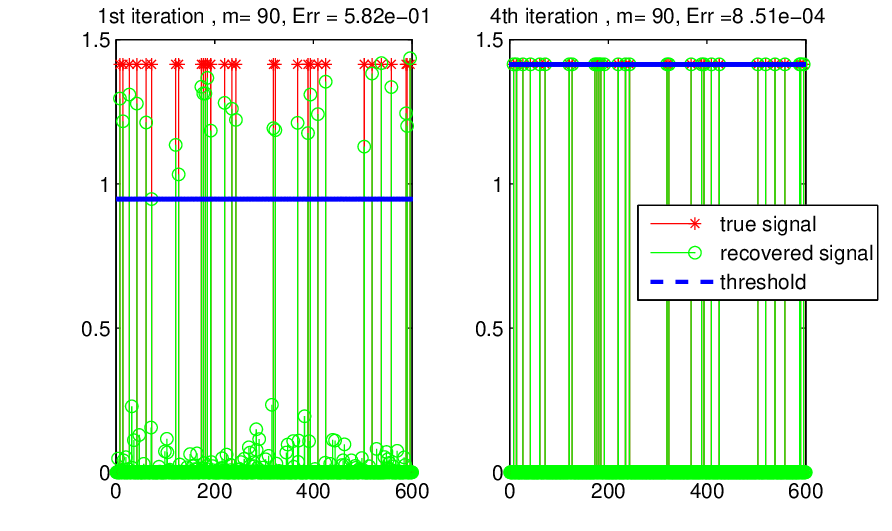}}
    \subfigure[]{
   \includegraphics[scale=0.43]{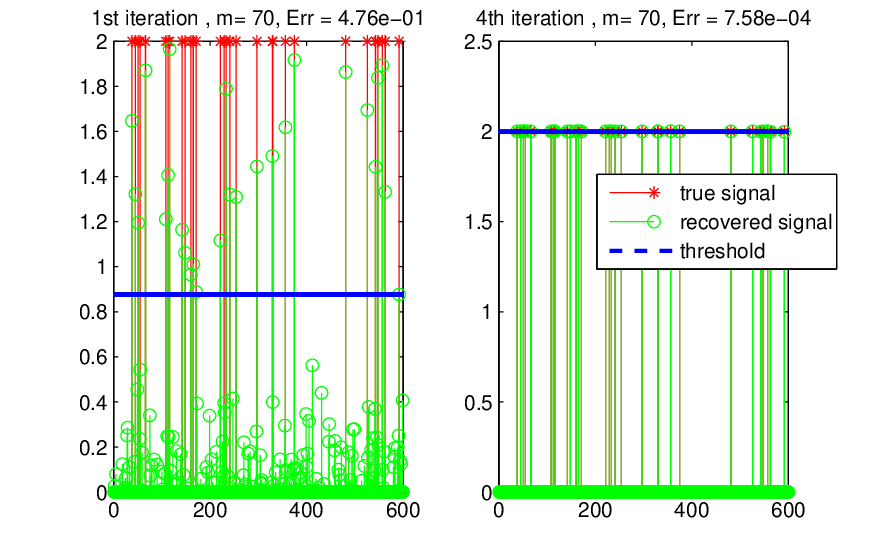}}
    \subfigure[]{
   \includegraphics[scale=0.43]{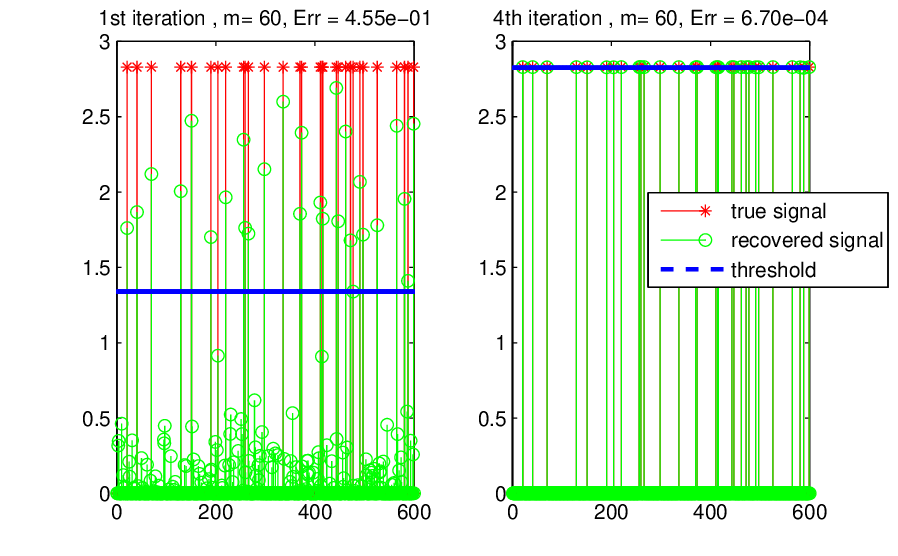}}
    \subfigure[]{
   \includegraphics[scale=0.43]{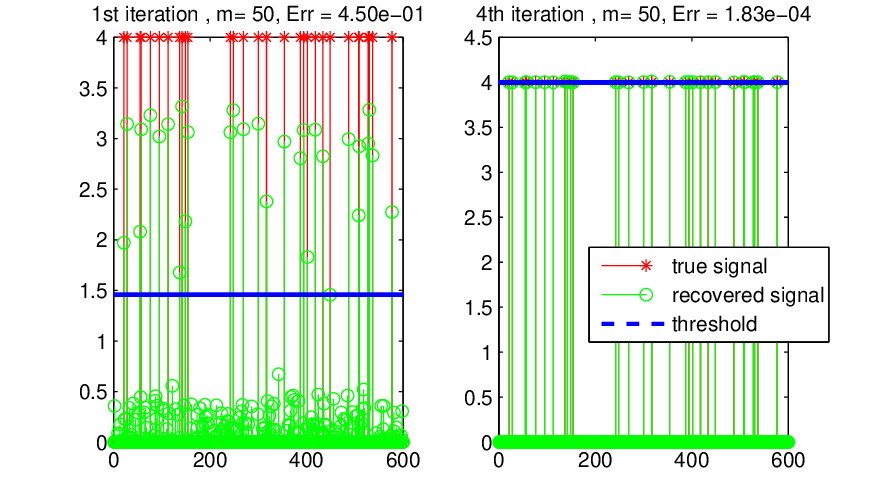}}
   \centering
   \begin{tabular}{|c|c|c|c|c|c|}
\hline 
 \multirow{2}{*}{L--Iteration} & \multicolumn{4}{|c|}{Nonzeros} &  \multirow{2}{*}{Relative error}     \\
 \cline{2-5}
     & Total true & Detected & Correct & False &  \\
    \cline{1-6}
     1--1 &  30 & 33	& 27   &  6  & 7.69e-01 \\

     1--4 &  30 & 36	& 30   &  6  & 6.46e-01 \\
     \cline{1-6}

     2--1 &  30 & 32	& 30  &  2  & 5.82e-01 \\

     2--4 &  30 & 30 & 30  &  0  & 8.51e-04 \\
    \cline{1-6}
     4--1 &  30 & 35	& 30  &  5  & 4.76e-01 \\

     4--4 &  30 & 30 & 30  &  0  & 7.58e-04 \\
    \cline{1-6}
     8--1 &  30 & 35	& 30  &  5  & 4.56e-01 \\

     8--4 &  30 & 30 & 30  &  0  & 6.70e-04 \\
    \cline{1-6}
     16--1 &  30 & 32  & 30 &  2  & 4.50e-01 \\

     16--4 &  30 & 30 & 30  &  0  & 1.83e-04 \\
    \cline{1-6}
\hline
\end{tabular}
   \caption{\small Compare the true Bernoulli signals and recovered signals obtained by ISDJS in different channels, where the two components in each subplot are the results in the first iteration and the fourth iteration respectively. (a)L=1, m=110, (b)L=2, m=90, (c)L=4, m=70, (d)L=8, m=60, (e)L=16, m=50.}
   \label{fig:6}
\end{figure}

{ We show a surprising performance of ISD with joint sparsity. In \cite{Wang10ISD} for the single channel signal recovery,
the $threshold$-ISD works well for signals with a fast decaying property of nonzero entries such as sparse Gaussian signals and certain power-law decaying signals. However, it does not work for signals that decay slowly or have no decay at all such as sparse Bernoulli signals, since the threshold based support detection fail to accurately distinguish true nonzero components according to the intermediate recovery results.}

{ As Fig \ref{fig:6} (a) presented, the support detection is poor and fails to correctly detect the true nonzero components in the single channel signal recovery. Nevertheless,  Fig \ref{fig:6} (b) shows the threshold based support detection can accurately find some true nonzero components,  even just for $L=2$. Then, in Figs \ref{fig:6} (c), (d) and (e), threshold based support detection works well. Finally, the ISDJS achieves quite good recovery performance, which suggests that ISD is able to achieve relatively accurate support detection even for Bernoulli signals by incorporating the joint sparsity structure, as briefly explained in Section \ref{sec:sd}. The table below Fig \ref{fig:6} displays that the support detection works well as the iteration proceeds. }

\begin{figure}[!hbtp]
  \centering
   \includegraphics[scale=0.44]{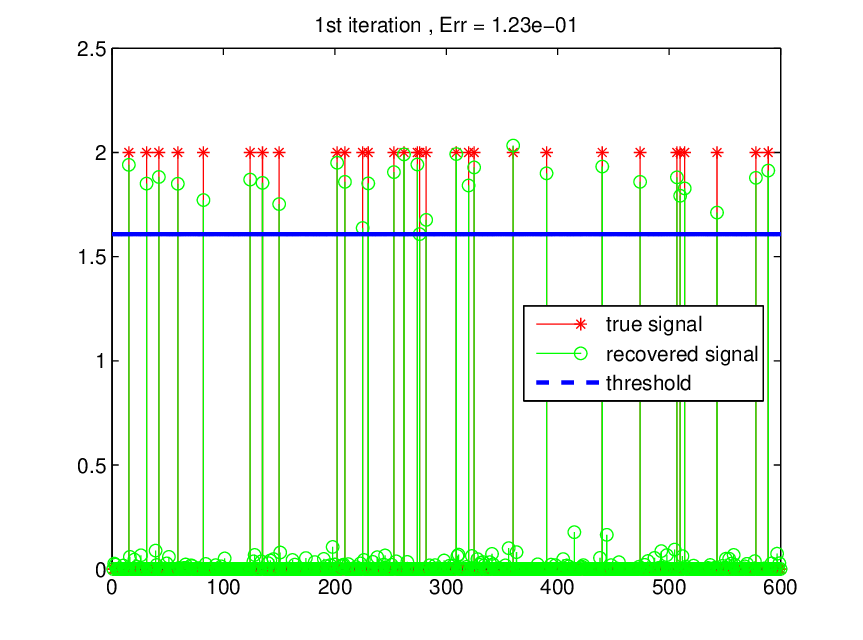}
   \includegraphics[scale=0.44]{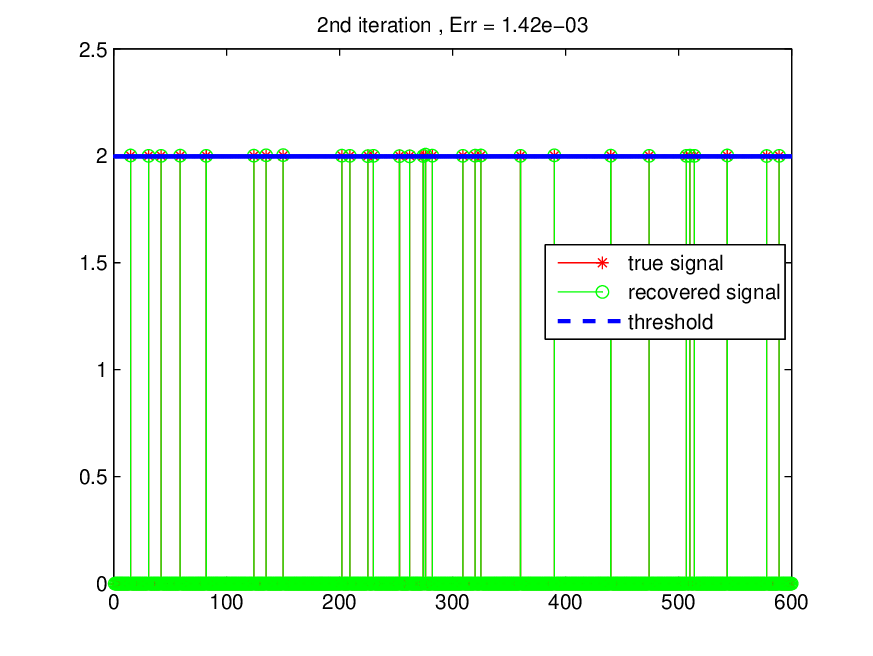}
   \includegraphics[scale=0.44]{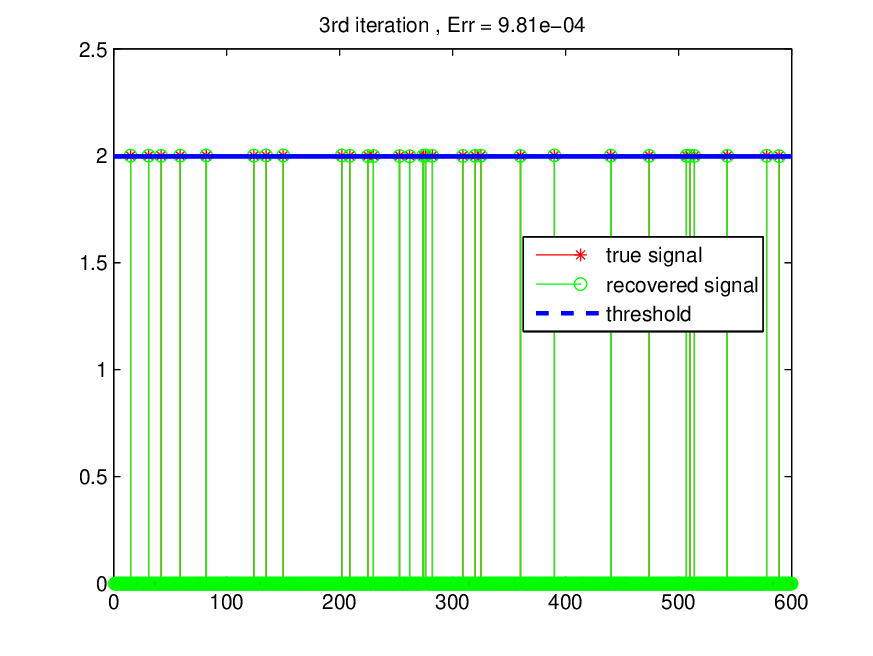}
   \includegraphics[scale=0.44]{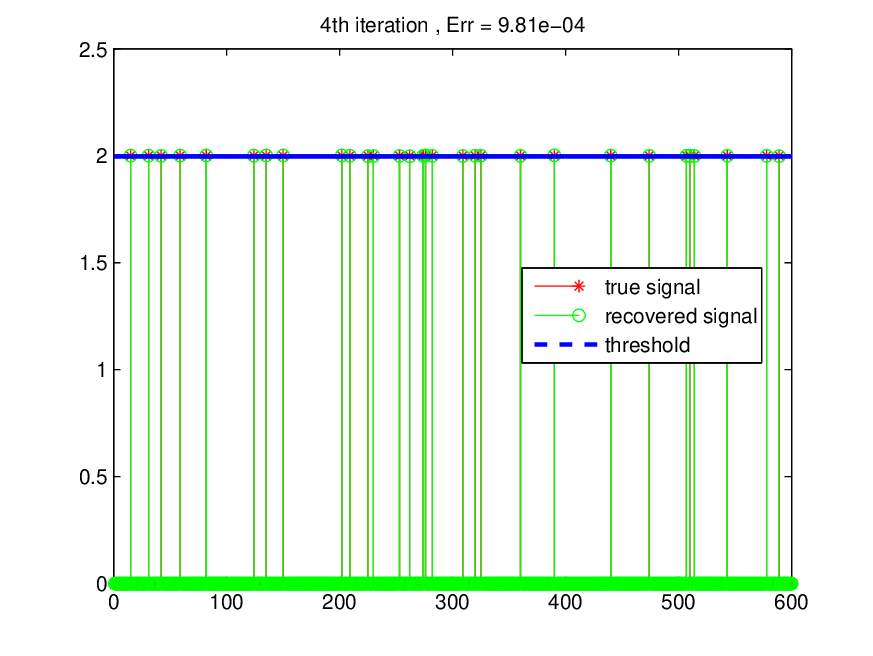}
\centering
\begin{tabular}{|c|c|c|c|c|c|}
\hline 
 \multirow{2}{*}{Iteration} & \multicolumn{4}{|c|}{Nonzeros} &  \multirow{2}{*}{Relative error}     \\
 \cline{2-5}
     & Total true & Detected & Correct & False &  \\
    \cline{1-6}
     1 &  30 & 28	& 26  &  2  & 1.23e-01 \\

     2 &  30 & 29	& 29  &  0  & 1.42e-03 \\

     3 &  30 & 30	& 30  &  0  & 9.81e-04 \\

     4 &  30 & 30   & 30  &  0  & 9.81e-04 \\
    \cline{1-6}
\hline
\end{tabular}
   \caption{\small Compare the true Bernoulli signals and recovered signals from ISDJS in each iteration with the L=4 channels.}
   \label{fig:7}
\end{figure}

\begin{figure}[!htbp]
  \centering
    \subfigure[]{
   \includegraphics[scale=0.43]{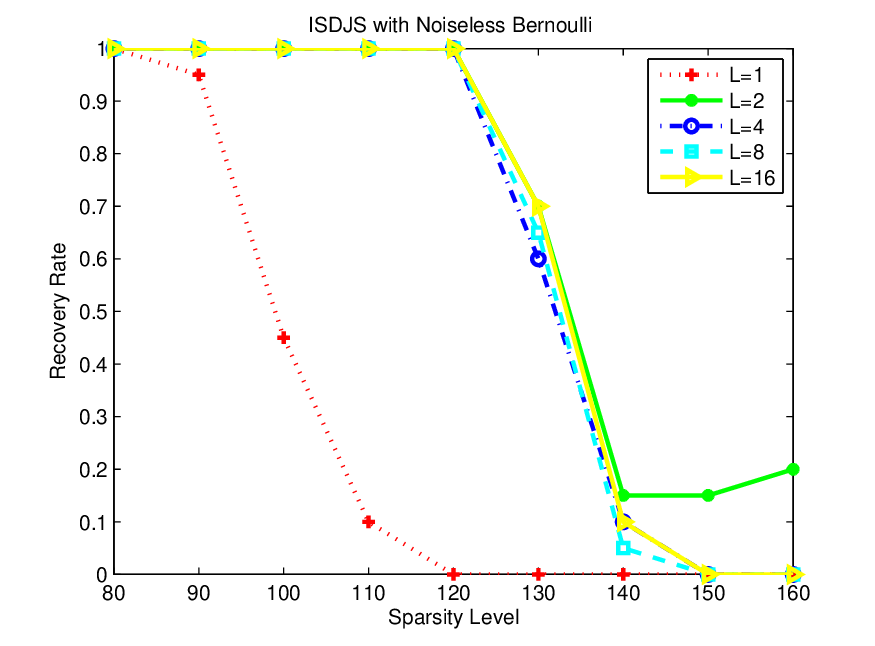}}
    \subfigure[]{
   \includegraphics[scale=0.43]{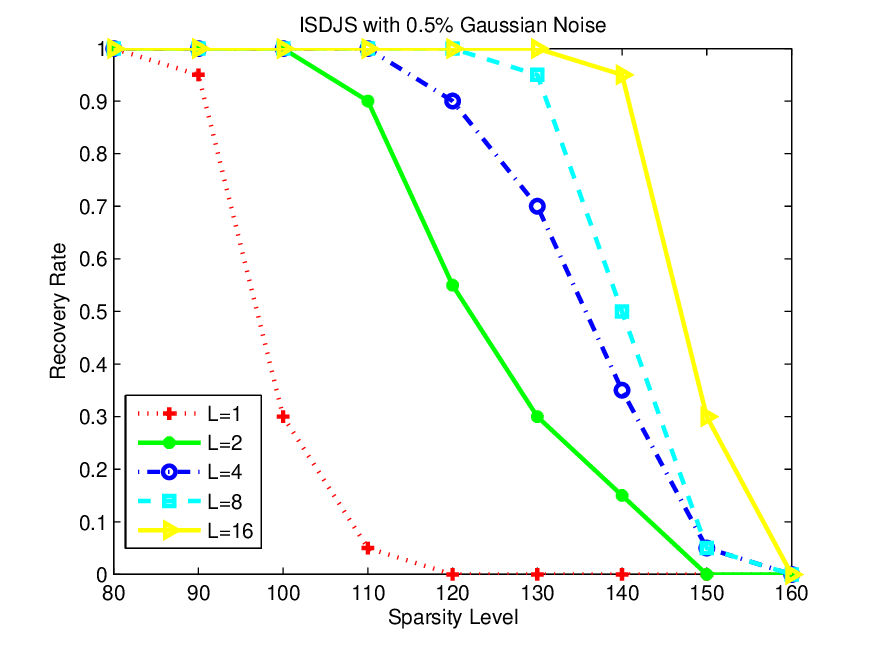}}
    \subfigure[]{
   \includegraphics[scale=0.43]{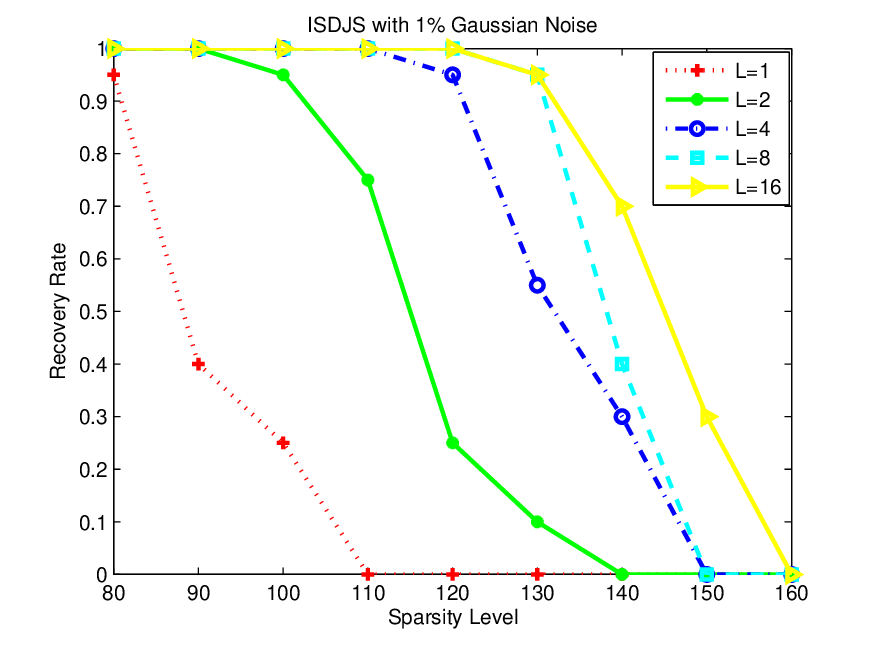}}
    \subfigure[]{
   \includegraphics[scale=0.43]{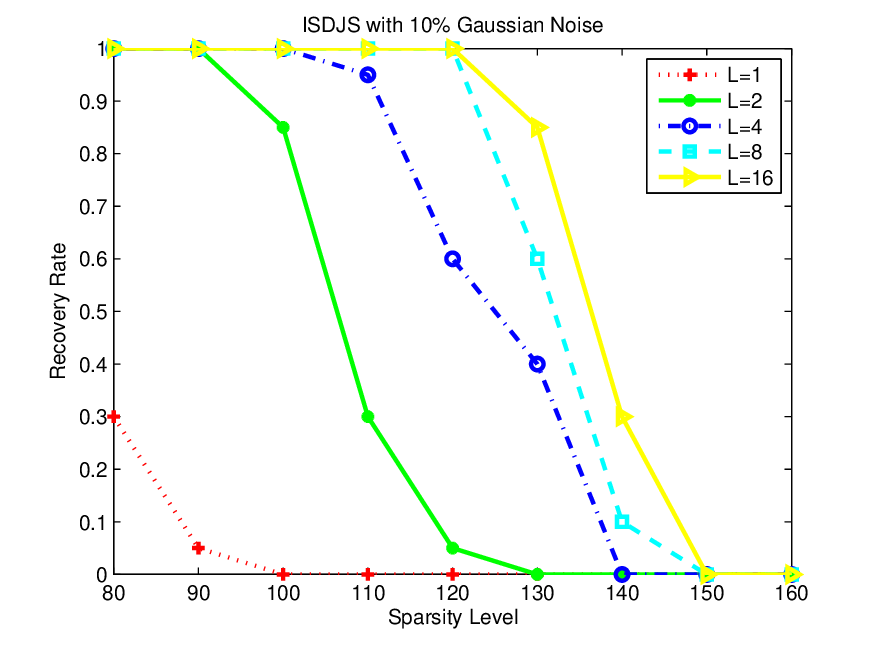}}
   \caption{\small Compare the recovery rate of ISDJS with L=1, 2, 4, 8, 16 in different noise levels for Bernoulli signals, (a)noiseless, (b)0.5\% noise, (c)1\% noise, (d)10\% noise.}
   \label{fig:8}
\end{figure}

In Fig \ref{fig:7}, we exhibit the performance of ISDJS for sparse Bernoulli signals of each outer iteration by taking $L=4$ and $m=70$ as an example.  It is possible to add more iterations but four iterations are  enough for ISDJS to return an accurate solution.
In Fig \ref{fig:8}, we show the recovery rates of  ISDJS with different channel number settings under four noise levels. The ISDJS consistently performs much better in multichannel cases rather than the single channel situation. Moreover, the ISDJS keeps robust in different noise levels.

\begin{figure}[!htbp]
  \centering
   \subfigure[]{
   \includegraphics[scale=0.33]{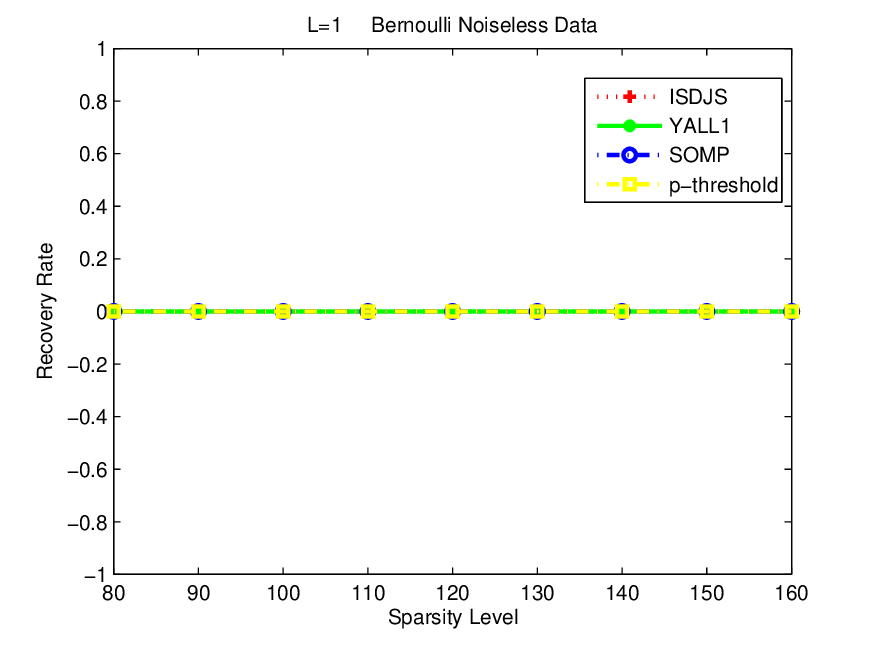}}
    \subfigure[]{
   \includegraphics[scale=0.33]{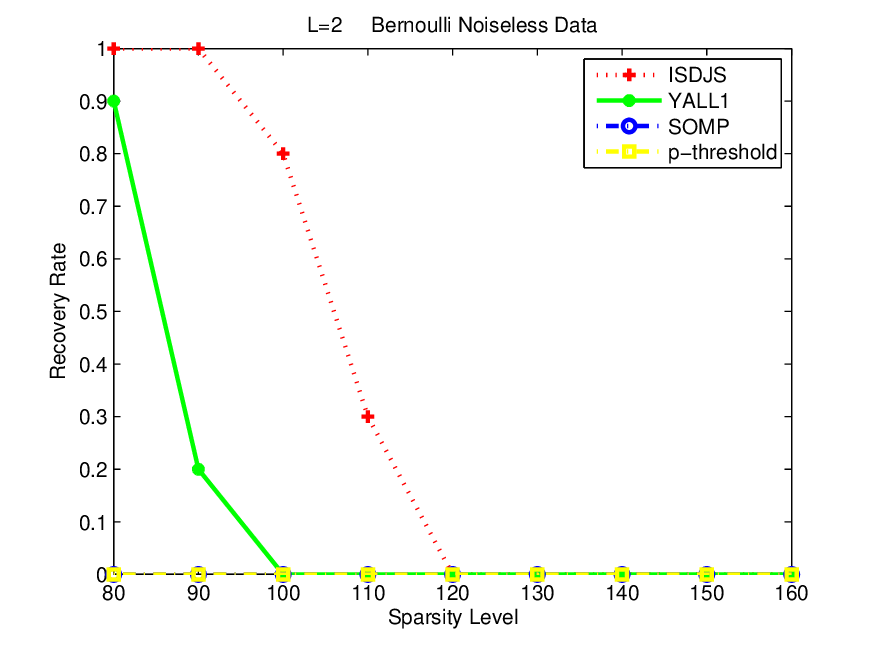}}\\
   \includegraphics[scale=0.305]{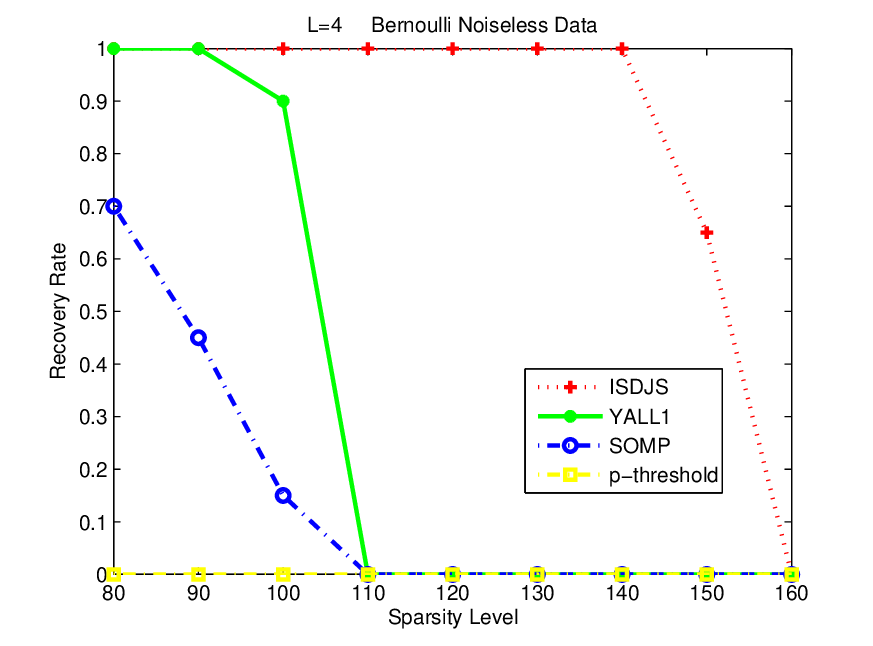}
   \includegraphics[scale=0.305]{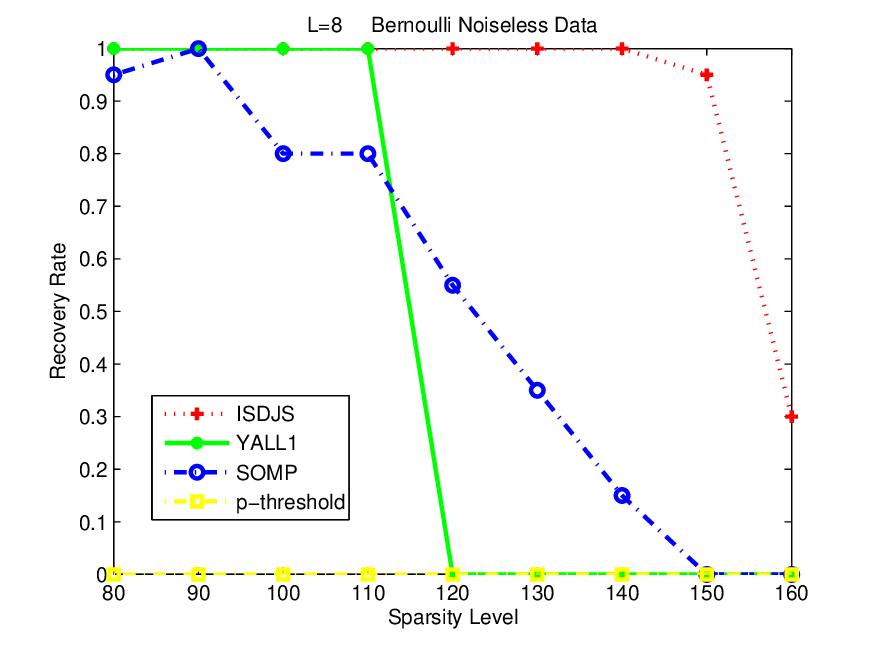}
   \includegraphics[scale=0.3]{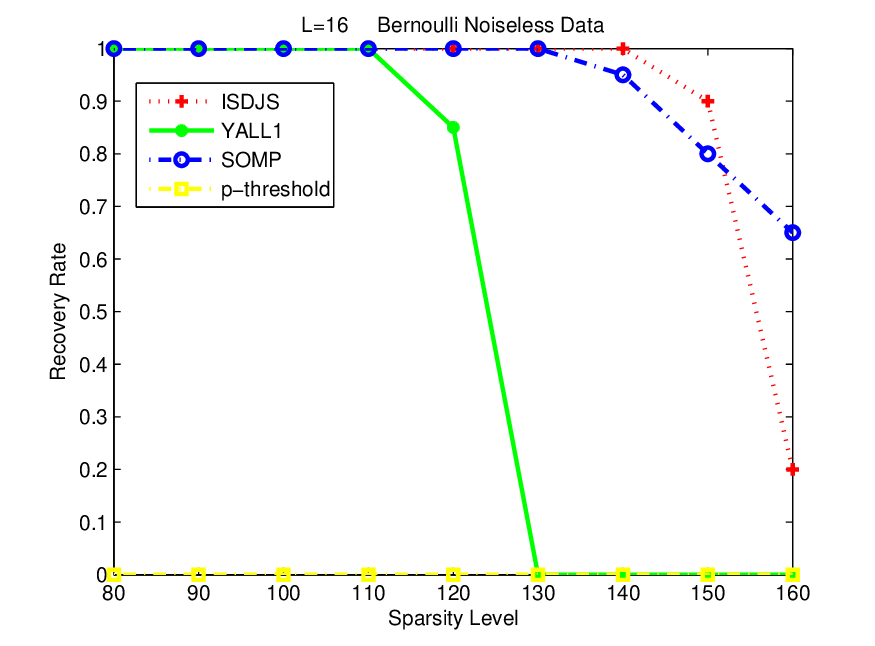}\\
   \hspace{0.05in}\small(c)\hspace{1.5in}(d)\hspace{1.5in}(e)
   \caption{\small Compare the recovery rate of four algorithms in different channels for noiseless Bernoulli signals, (a)L=1, (b)L=2, (c)L=4, (d)L=8, (e)L=16.}
   \label{fig:9}
\end{figure}

\begin{figure}[!htbp]
  \centering
    \subfigure[]{
   \includegraphics[scale=0.32]{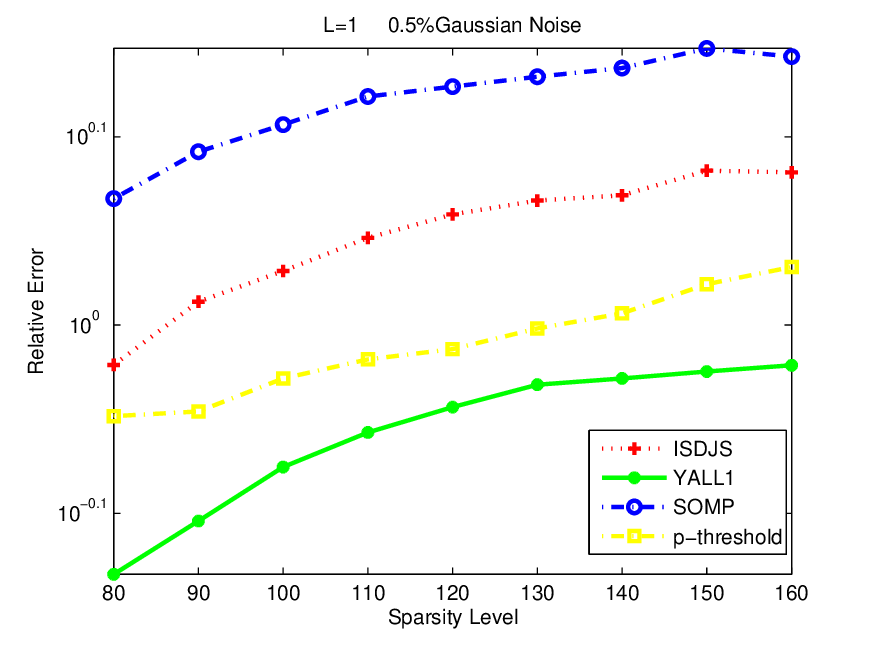}}
    \subfigure[]{
   \includegraphics[scale=0.32]{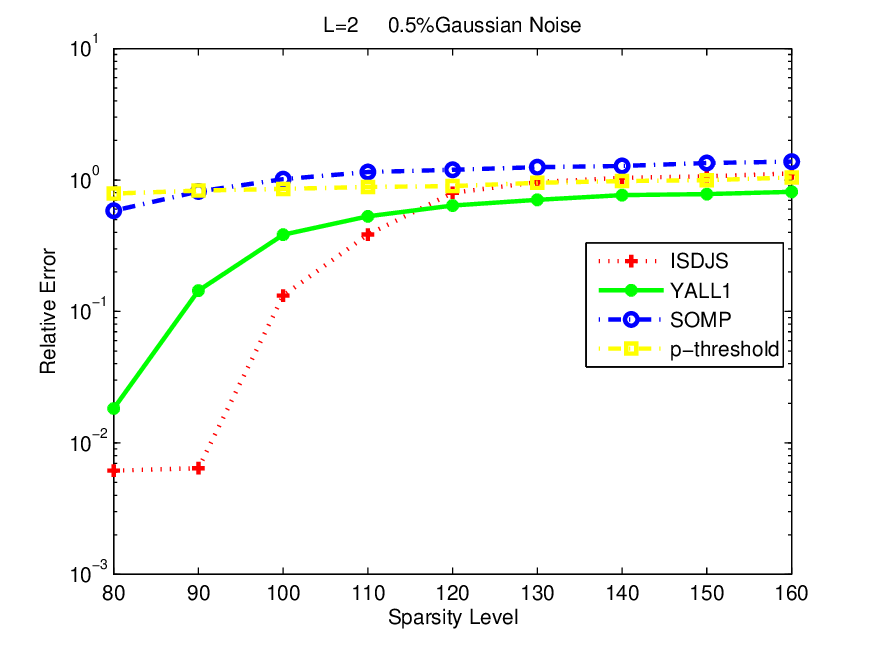}}\\

   \includegraphics[scale=0.305]{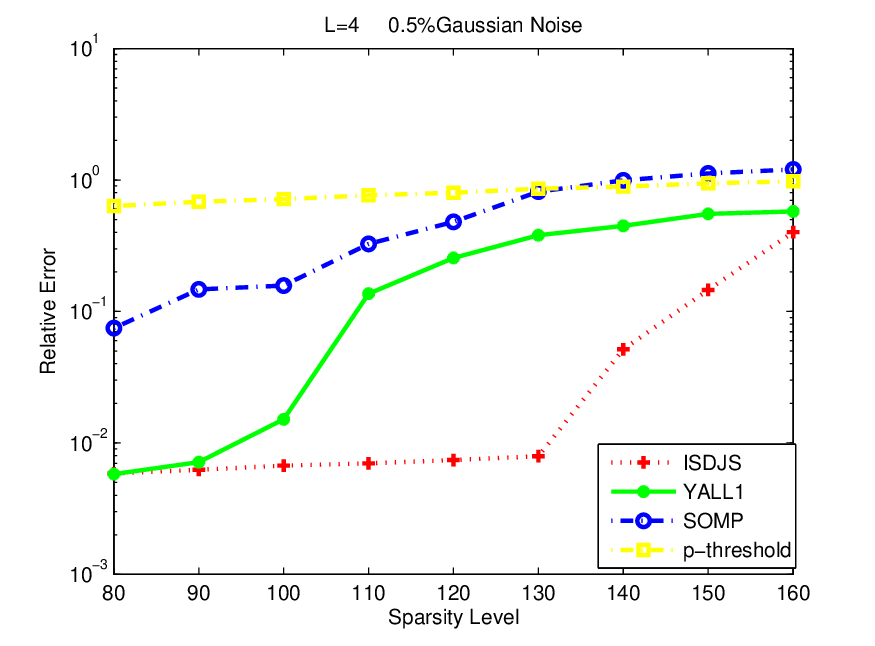}
   \includegraphics[scale=0.305]{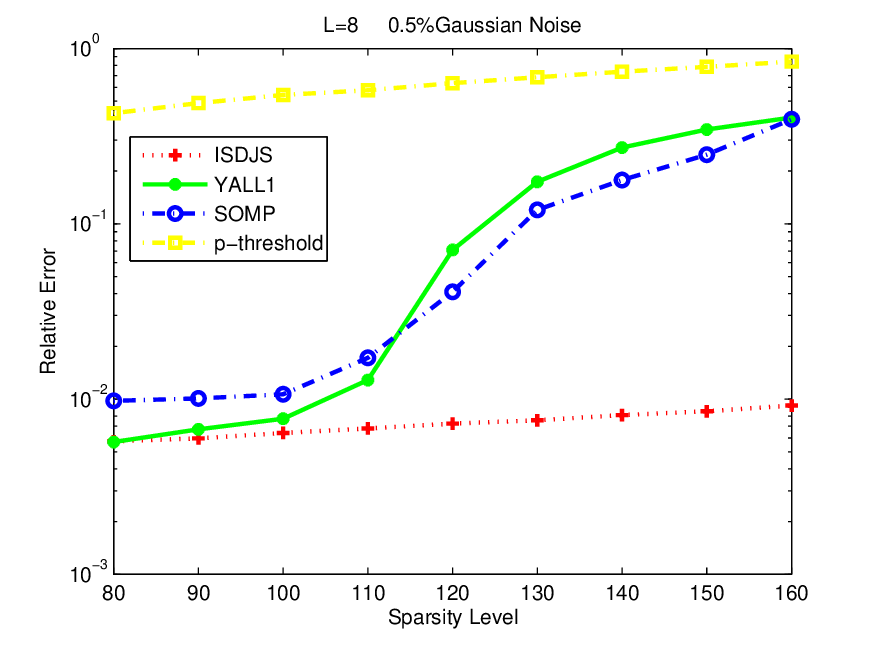}
   \includegraphics[scale=0.3]{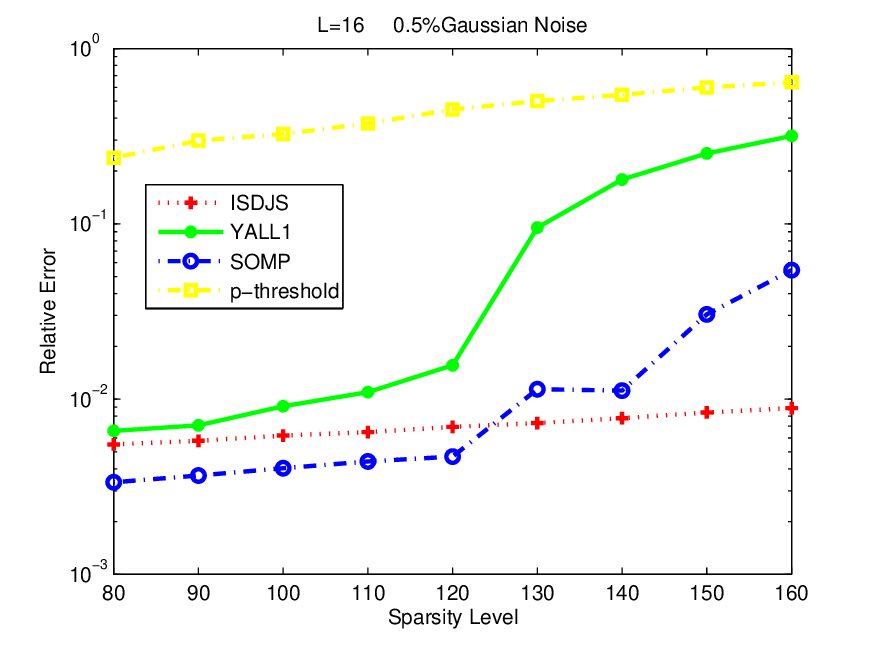}\\
   \hspace{0.05in}\small(c)\hspace{1.5in}(d)\hspace{1.5in}(e)
   \caption{\small Compare relative error of four algorithms in different channels with $0.5\%$ noise for Bernoulli signals, (a)L=1, (b)L=2, (c)L=4, (d)L=8, (e)L=16.}
   \label{fig:10}
\end{figure}

We plot the recovery rate of ISDJS in comparison with other three algorithms for noiseless sparse Bernoulli signals in Fig \ref{fig:9}.
Obviously, Fig \ref{fig:9} (a) shows that all test algorithms perform poor on single channel Bernoulli signals, since  the joint structure prior of signals do not exist here. Surprisedly, the recoverability of ISDJS is dramatically improved as the channel numbers increases in Figs \ref{fig:9} (b), (c), (d) and (e). Similarly, Fig \ref{fig:10} exhibits the relative error of all test algorithms with $0.5\%$ noise in different channel number settings.

All above numerical experiments  attest that ISDJS can make significant improvement for multichannel sparse signal recovery even without the fast decaying property, by incorporating joint sparsity property  into the implementation of threshold based support detection.
\subsection{An Example from Collaborative Spectrum Sensing}
Now we consider a compressive spectrum sensing scheme for cognitive radio networks \cite{J2011,chang2015}.  Spectrum sensing aim to detect spectrum holes (i.e., channels not used by any primary users). The cognitive radio (CR) nodes must constantly sense the spectrum in order to detect the presence of the primary radio (PR) nodes and use the spectrum holes without causing harmful interference to the PRs. In practice, improving the ability of detecting complete spectrum usage in collaborative spectrum sensing is an important topic  but also a major challenge in CR networks.

\begin{figure}[!htbp]
  \centering
   \subfigure[]{
   \includegraphics[scale=0.4]{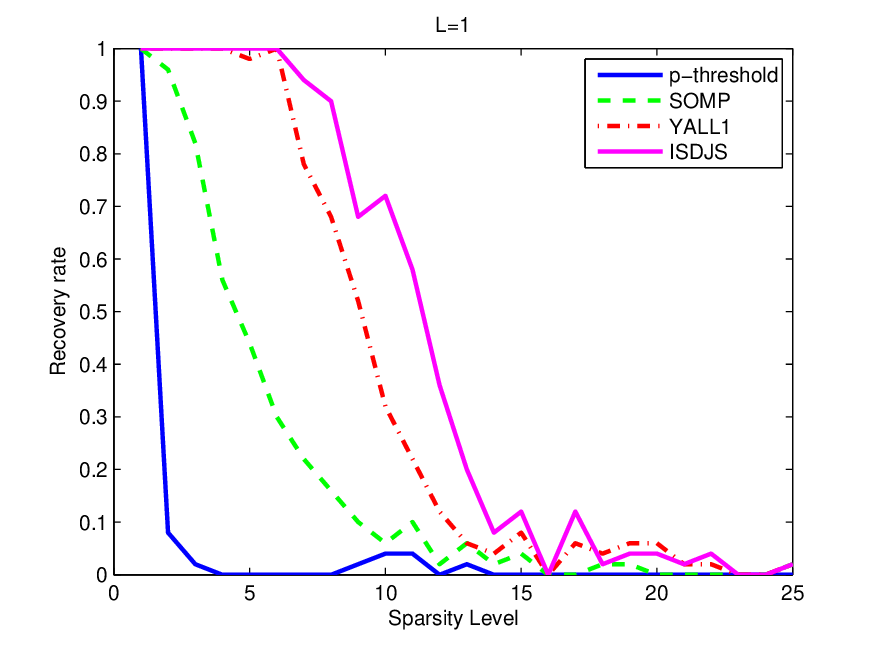}}
   \subfigure[]{
   \includegraphics[scale=0.4]{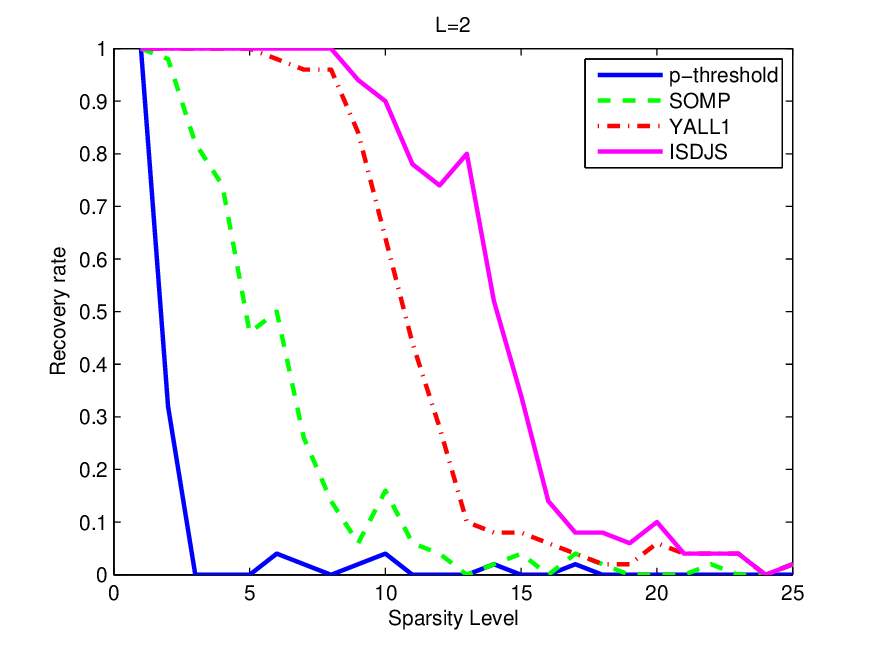}}
    \subfigure[]{
   \includegraphics[scale=0.4]{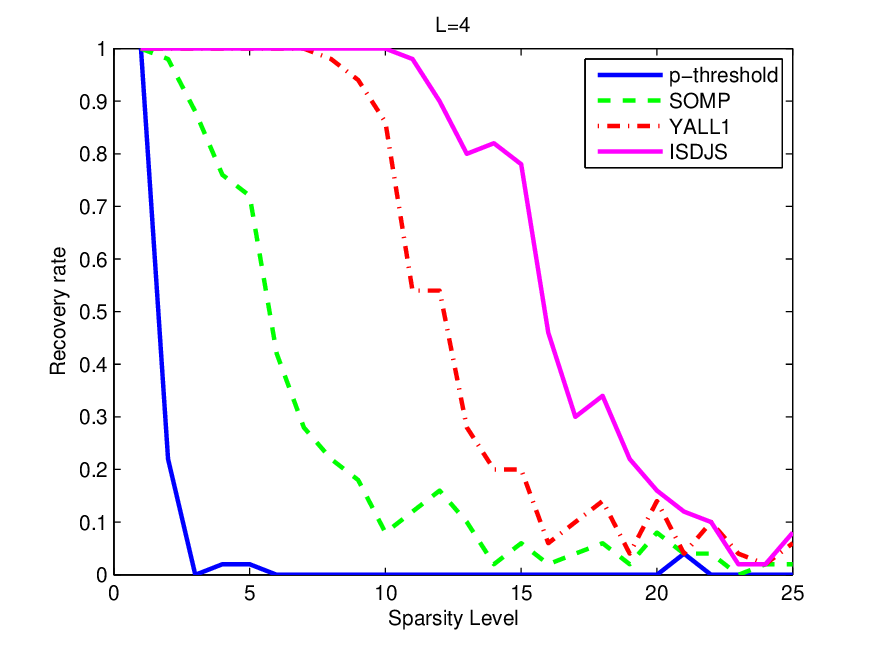}}
    \subfigure[]{
   \includegraphics[scale=0.4]{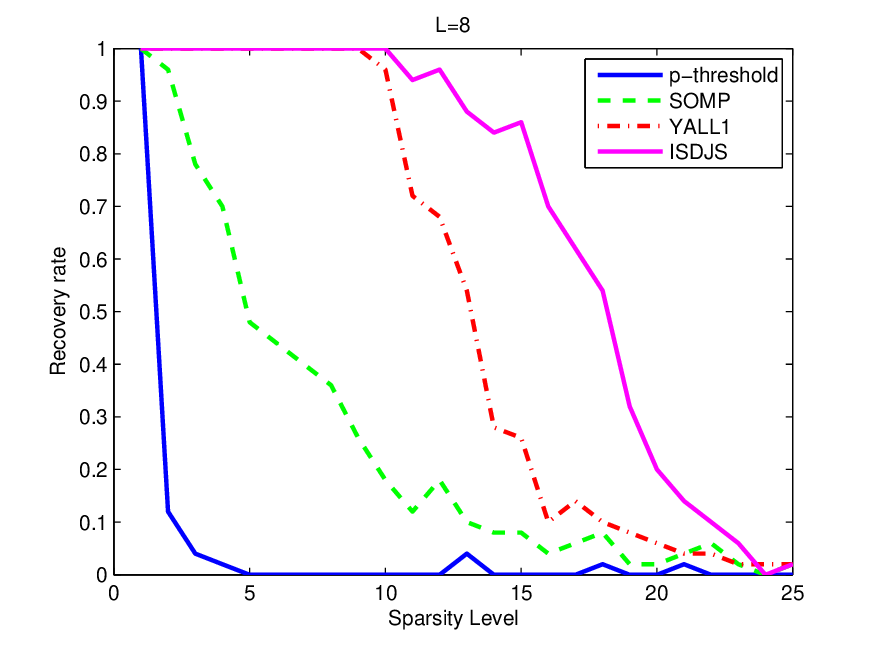}}
    \subfigure[]{
   \includegraphics[scale=0.4]{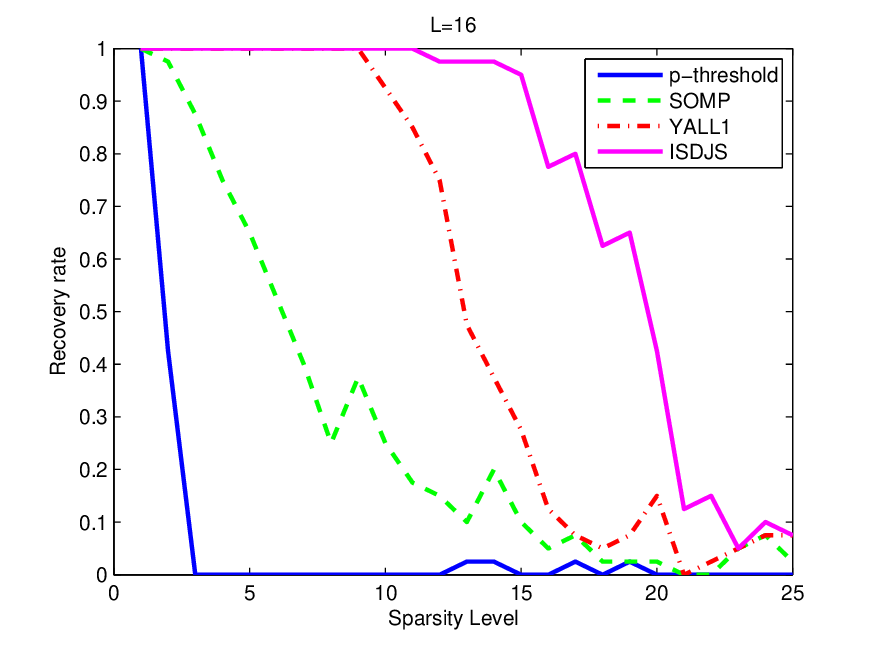}}
   \caption{\small Compare the recovery rate of four algorithms in different channel numbers in spectrum sensing, (a)L=1, (b)L=2, (c)L=4, (d)L=8, (e)L=16.}\label{fig:UH}
\end{figure}

We view a $l$-node cognitive radio network within a $500 \times 500$ meter square area centered at the fusion center. The $l$ CR nodes are uniformly randomly located. These cognitive radio nodes collaboratively sense the existence of primary users within a $1000 \times 1000$ meter square area on $n$ channels, which are centered also at the fusion center. A channel is either occupied by a PR or unoccupied, corresponding to the states 1 and 0, respectively.
Let an $n \times n$ diagonal matrix $H$ represent the states of all the channel sources using 0 and 1 as diagonal entries, indicating the unoccupied or occupied states, respectively.  Channel gains are characterized by an $l \times n$ channel gain matrix $G$. Then, the collaborative spectrum sensing model can be formulated as follows \cite{J2011}:
\begin{equation}
X_{n\times l} = H_{n\times n}(G_{l\times n})^{T}.
\end{equation}
For $X$, the $j$-th column of $X$ corresponds to the channel occupancy status received by the $j$-th CR, and the $i$-th row of $X$ corresponds to the occupancy status of the $i$-th channel. A row has a positive value if and only if the $i$-th channel is used. Since there are only a small number of used channels, $X$ is sparse in terms of the number of nonzero rows.
In this example, we set $n=25$ and $l=1, 2, 4, 8, 16$. We apply the ISDJS to solve above collaborative spectrum sensing model.
Fig \ref{fig:UH} presents the results of ISDJS compared with YALL1 group, SOMP and p-threshold algorithms in different settings of $l$. With the sparsity level (nonzero rows) of $X$ increasing, the advantage of ISDJS is notable.

\subsection{An Example for Multi-task Feature Learning}\label{sec:multilearning}

We provide an example to demonstrate the performance of ISDJS for multi-task feature learning. A real-world data set, i.e. Letter \cite{JLiu,G07} is employed in this experiment.  The Letter data set was collected by the MIT Spoken Language Systems Group \footnote{http:www.seas.upenn.edu/~taskar/ocr/}. It contains 8 default tasks for the handwritten letters. The writings are collected from over 180 different writers and there are 45,679 samples, where the letters are represented by $8\times16$ binary pixel images. This is a typical multi-task feature learning problem.

\begin{figure}[!htbp]
  \centering
   \includegraphics[scale=0.46]{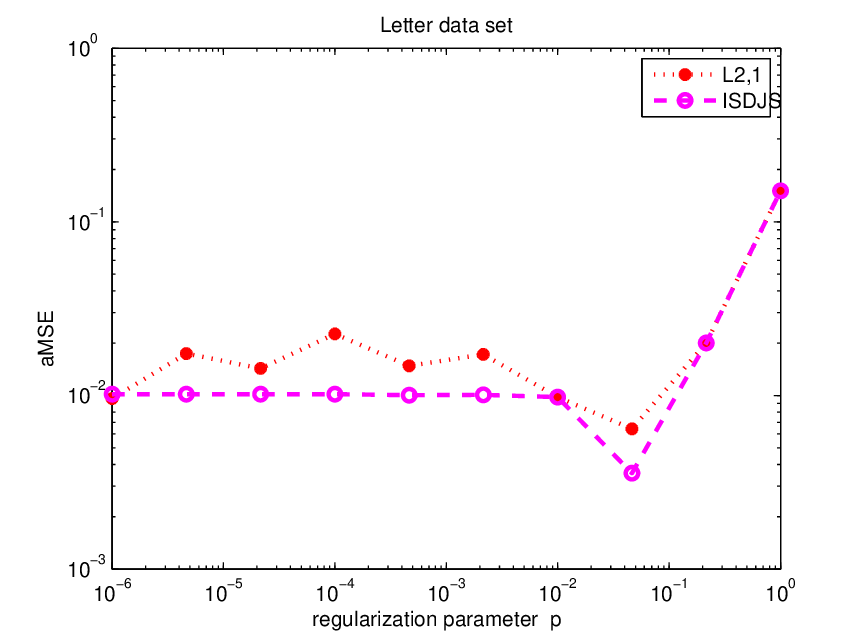}
   \includegraphics[scale=0.45]{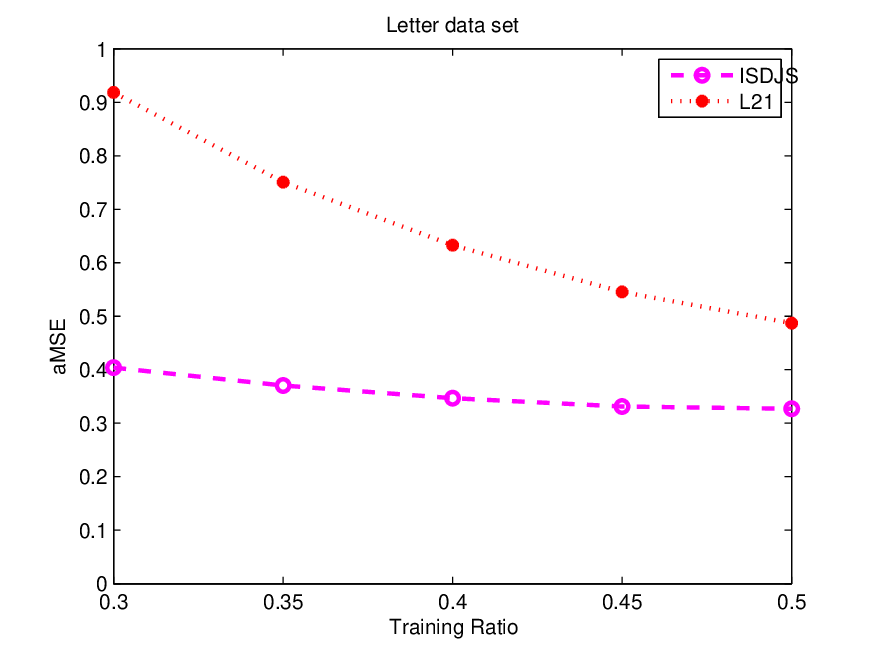}\\
   \hspace{0.7in}\small(a)\hspace{2.5in}(b)\hspace{0.7in}
   \caption{ Compare the algorithm in \cite{JLiu} with ISDJS on Letter data set. (a) aMSE vs. regularization parameter $\rho$, (b) aMSE vs. training ratio.}
   \label{fig:multitask}
\end{figure}

The baseline algorithm for solving the common $\ell_{2,1}$ regularized model proposed in  \cite{JLiu} is used to compare with ISDJS. Here we do not evaluate the performance according to the estimation error of the weight matrix $X$, whose true values are unknown in practice. Instead, we use the averaged means squared error (aMSE) and normalized mean squared error (nMSE):
$$
\mathrm{aMSE}=\frac{||\hat{b}-\bar{b}||_{F}}{||\bar{b}||_{F}},
$$
$$
\mathrm{nMSE}=N\frac{||\hat{b}-\bar{b}||_{F}^{2}}{||\hat{b}||_{1}\cdot||\bar{b}||_{1}},
$$
where $\hat{b}$ is the predictive value of the trained model for the test set, $\bar{b}$ is the known reference true value and $N$ is the number of the test sample. Both nMSE and aMSE are commonly used in multi-task learning problems \cite{MSMTFL}.
\begin{figure}[!htbp]
  \centering
   \includegraphics[scale=0.45]{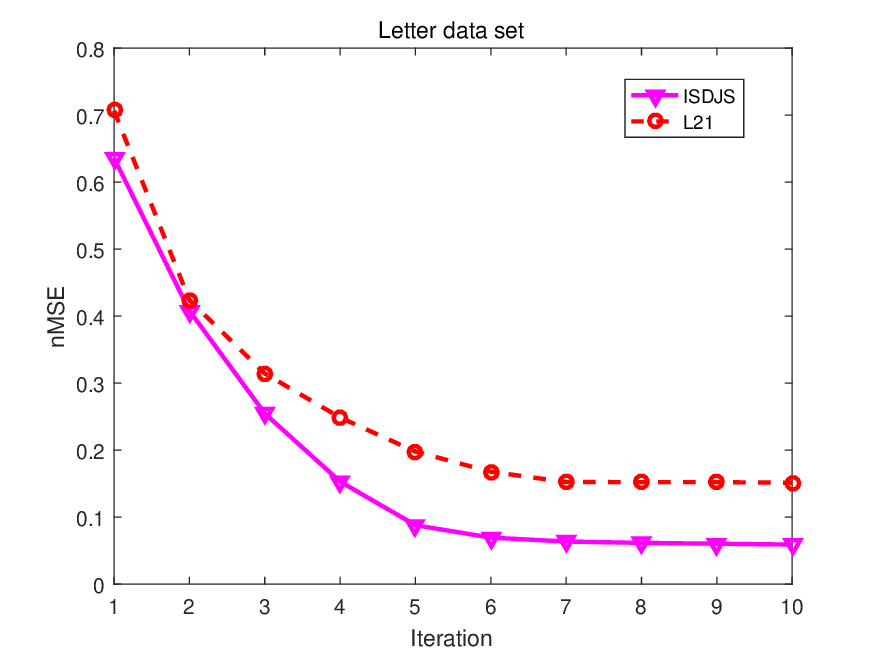}
   \includegraphics[scale=0.45]{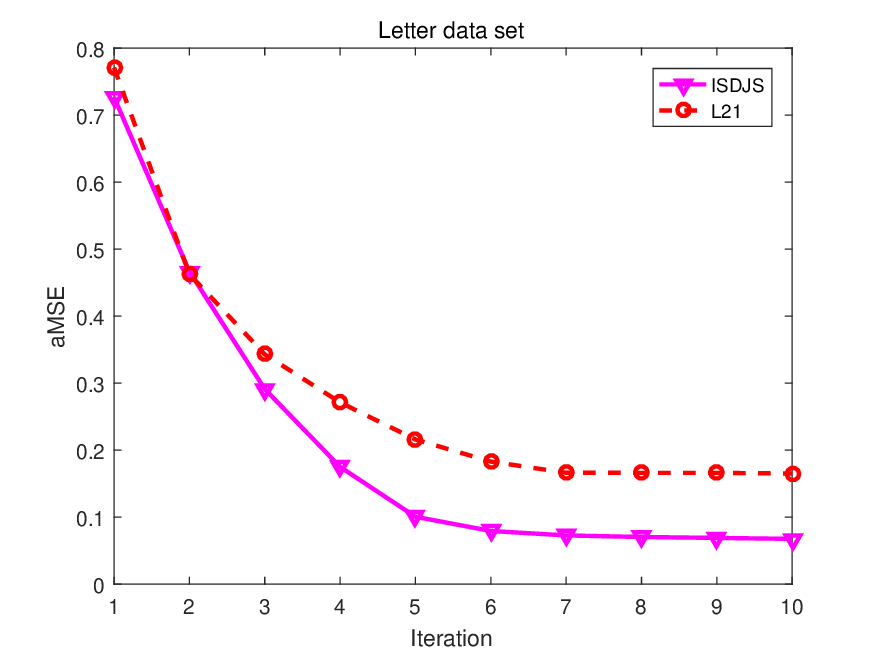}\\
   \hspace{0.7in}\small(a)\hspace{2.5in}(b)\hspace{0.7in}
   \caption{\small Compare the algorithm in \cite{JLiu} with ISDJS on Letter data set. (a) nMSE vs. iteration, (b) aMSE vs. iteration.}
   \label{fig:multi}
\end{figure}

It is well known that an appropriate regularization parameter is vital for the great performance of the algorithm.  As depicted in Fig \ref{fig:multitask} (a), ISDJS achieves a small error when the regularization parameter $\rho$ is around $4.64\times10^{-2}$. Then, we randomly extract the training samples from each task with different training ratios ($30\%$, $35\%$, $40\%$, $45\%$ and $50\%$) and exploit the rest of samples to form a test set. In Fig \ref{fig:multitask} (b), ISDJS performs much better even with a small training ratio. It is easy to observe that ISDJS is convergent after a few iterations in Fig \ref{fig:multi}, which is consistent with the $\mathbf{Theorem}$ $\mathbf{1}$ in Section \ref{sec:theory}.
This real-world experimental results for the multi-task feature learning further support the effectiveness of ISDJS.

\section{Conclusion} \label{sec:conclusion}
In this paper, we have proposed a truncated joint sparsity model and developed an efficient algorithm named ISDJS to enhance the sparse estimation performance. They are applied in the fields of compressive sensing and feature learning. The proposed method is an extension of self-learning based iterative support detection (ISD) from common sparsity to joint sparsity. The joint sparsity structure is naturally incorporated into the implementation of threshold based support detection and in this way the fast decaying property is no longer required. Then, we have elaborated some preliminary results of the convergence analysis and a sufficient recovery condition for the proposed method. Both synthetic and practical experiments demonstrate the better performance of ISDJS compared with several state-of-the-art alternatives. In the future, we will explore more applications such as image inpainting, image classification and feature selection to employ the ISDJS method, and design specific implementations of support detection to achieve the outstanding performance for different applications.

\section*{Acknowledgment}
The authors would like to thank Prof. Wenxing Zhang for valuable discussion about the theoretical analysis of the proposed method, and the editor and referees for their valuable suggestions and comments.
This work is supported by the 973 project (No. 2015CB856000), the Natural Science Foundation of China (91330201) and the Fundamental Research Funds for the Central Universities (ZYGX2013Z005).


\section*{References}

\bibliographystyle{plain}
\bibliography{EJS}

\end{document}